\documentclass[12pt]{article}
\usepackage{srcltx}
\usepackage{amssymb}
\usepackage{amsfonts}
\usepackage{amsmath}
\usepackage{amsthm}
\usepackage{enumerate}
\usepackage{multicol}

\newtheorem{theorem}{Theorem}[section]
\newtheorem{prop}[theorem]{ Proposition}
\newtheorem{lemma}[theorem]{Lemma}

\newtheorem{remark}{Remark}[section]

\newcommand\cA{{\cal A}}
\newcommand\cC{{\cal C}}
\newcommand\cE{{\cal E}}

\newcommand\cL{{\cal L}}
\newcommand\cB{{\cal B}}

\newcommand\cN{{\cal N}}

\newcommand\cP{{\cal P}}

\newcommand\cR{{\cal R}}
\newcommand\cT{{\cal T}}

\newcommand\e{\epsilon}
\newcommand\ve{\varepsilon}

\newcommand\ov{\overline}

\def\bbr{{\mathbb R}}

\def\text#1{\hbox{#1}}

\def\endproof{\mbox{\ $\qed$}}

\def\E{{\bf E}}
\def\e{{\bf e}}
\def\A{{\bf A}}
\def\P{{\bf P}}

\def\H{{\bf H}}
\def\L{{\bf L}}

\newcommand{\wh}{\widehat}
\newcommand{\wt}{\widetilde}

\def\O{\hbox{\rm O}}

\def\Chi{{\bf 1}}
\def\d{\mathrm{d}}
\def\build #1_#2{\mathrel{\mathop{\kern 0pt #1}\limits_{#2}}} 
\newcommand{\zs}[1]{{\mathchoice{#1}{#1}{\lower.25ex\hbox{$\scriptstyle#1$}}
{\lower0.25ex\hbox{$\scriptscriptstyle#1$}}}}
\numberwithin{equation}{section}

\begin{document}
\title{Adaptive asymptotically efficient estimation\\
in heteroscedastic nonparametric regression via model selection.
}
\author{{\Large By  Leonid Galtchouk and Sergey Pergamenshchikov}\\
Louis Pasteur University  of Strasbourg and University of Rouen
}
\maketitle

\begin{abstract}
The paper deals with asymptotic properties of the adaptive
procedure  proposed in the author paper, 2007, 
for  estimating a unknown nonparametric regression.
We prove that this
procedure is asymptotically efficient for a quadratic risk,
i.e. the asymptotic quadratic risk for this procedure
coincides with the Pinsker constant which gives a sharp 
 lower bound for the quadratic risk over all possible estimates.
\end{abstract}
\bibliographystyle{plain}
\renewcommand{\columnseprule}{.1pt}
\newpage

\section{Introduction}\label{I}

The paper deals with the  estimation problem in 
the heteroscedastic nonparametic regression model 
\begin{equation}\label{I.1}
y_\zs{j}=S(x_\zs{j})+\sigma_\zs{j}(S)\,\xi_\zs{j}\,, 
\end{equation}
where the design points $x_\zs{j}=j/n$, $S(\cdot)$
 is an unknown function to be estimated, 
$(\xi_\zs{j})_\zs{1\le j\le n}$ is a sequence of centered 
 independent random variables with
unit variance and 
$(\sigma_\zs{j}(S))_\zs{1\le j\le n}$ are unknown scale functionals depending on
the design points and the regression function $S$.


 Typically, the notion of asymptotic
optimality is associated with the optimal convergence rate of the minimax risk
(see e.g., Ibragimov, Hasminskii,1981; Stone,1982).
  An important question in optimality
results is to study the exact asymptotic behavior of the minimax risk. Such results
have been obtained only in a limited number of investigations. As to the
nonparametric estimation problem
for heteroscedastic regression models we should mention the papers by Efromovich, 2007,
Efromovich, Pinsker, 1996, and Galtchouk, Pergamenshchikov, 2005,
 concerning the exact asymptotic behavior of the $\cL_\zs{2}$-risk
and the paper by Brua, 2007,
devoted to the efficient pointwise estimation for
heteroscedastic regressions. 

Heteroscedastic regression models are largely used in financial mathematics, in particular,
 in problem of calibrating (see e.g., Belomestny, Reiss, 2006). An example of 
heteroscedastic regression models is given by econometrics
(see, for example, Goldfeld, Quandt, 1972,
 p. 83), where 
for  consumer budget problems one  uses
 some parametric version of 
model \eqref{I.1}  with the scale coefficients  defined as
\begin{equation}\label{I.2}
\sigma^2_\zs{j}(S)=c_\zs{0}+c_\zs{1}x_\zs{j}+c_\zs{2}S^2(x_\zs{j})\,,
\end{equation}
where $c_\zs{0}$, $c_\zs{1}$ and $c_\zs{2}$ are some unknown positive  constants.

The purpose of the article is to study asymptotic properties of the adaptive 
estimation procedure proposed in Galtchouk, Pergamenshchikov, 2007,
for which a non-asymptotic oracle
inequality was proved for quadratic risks. We will prove that this oracle inequality
is  asymptotically sharp, i.e. the asymptotic quadratic risk is minimal. It means
the adaptive estimation procedure is efficient under some the conditions on the scales
$(\sigma_\zs{j}(S))_\zs{1\le j\le n}$ which are satisfied in the case \eqref{I.2}.
Note that in  Efromovich, 2007, Efromovich, Pinsker, 1996,
an efficient adaptive procedure is constructed
for heteroscedastic regression when the scale coefficient is independent of $S$, i.e.
$\sigma_\zs{j}(S)=\sigma_\zs{j}$. In Galtchouk, Pergamenshchikov, 2005,
for the model \eqref{I.1} the asymptotic efficiency
was proved under strong the conditions on the scales which are not satisfied in the case \eqref{I.2}.
Moreover in the cited papers the efficiency was proved for the gaussian random variables
$(\xi_\zs{j})_\zs{1\le j\le n}$ that is very restrictive for applications of proposed methods
to practical problems.

In the paper we modify the risk. We take a additional supremum over the family
of unknown noise distributions like to Galtchouk, Pergamenshchikov, 2006.
 This modification allows us to eliminate from
the risk dependence on the noise distribution. Moreover for this risk a efficient procedure 
is robust with respect to changing the noise distribution. 

It is well known to prove the asymptotic efficiency one has to show that the asymptotic quadratic 
risk coincides with the lower bound which is equal to the Pinsker constant. In the paper two 
problems
are resolved: in the first one a upper bound for the risk is obtained by making use of the
non-asymptotic oracle inequality from Galtchouk, Pergamenshchikov, 2007,
in the second one we prove that this upper
bound coincides with the Pinsker constant. Let us remember that the adaptive procedure proposed in
Galtchouk, Pergamenshchikov, 2007, is based on weighted least-squares estimates, where the 
weights are proper
modifications of the Pinsker weights for the homogeneous case (when
$\sigma_1(S)=\ldots=\sigma_\zs{n}(S)=1$) relative to a certain smoothness of the function $S$ and
this procedure chooses a best estimator for the quadratic risk among these estimators. To obtain 
the Pinsker constant for the model \eqref{I.1} one has to prove a sharp asymptotic lower bound 
for the quadratic risk in the case when the noise variance depends on the unknown regression
function.
In this case, as usually, we minorize the minimax risk by a bayesian one for a respective parametric family. Then for the bayesian risk we make use of a lower bound (see Theorem 6.1)
which is a modification of the van Trees inequality (see, Gill, Levit, 1995).

    The paper is organized as follows.  
In Section~\ref{sec:Ad} we construct an adaptive estimation procedure.
In Section~\ref{sec:Co} we formulate principal the conditions. The main results
are presented in Section~\ref{sec:Ma}. The upper bound for the quadratic risk is given
in Section~\ref{sec:Up}. 
In Section~\ref{sec:Lo} we give all main steps of proving the lower bound.
In Subsection~\ref{subsec:Tr} we find the lower bound for the bayesian
 risk which minorizes the minimax risk.
 In Subsection~\ref{subsec:Fa} we study a special parametric functions family 
used to define the bayesian risk.
In Subsection~\ref{subsec:Br} we choose a prior distribution for bayesian risk
to maximize the lower bound.
Section~\ref{sec:Np} is devoted to explain how to use the given procedure
 in the case when the unknown regression function is non periodic.
 In Section~\ref{sec:Cn} we discuss the main results and their practical importance.
The proofs are given in Section~\ref{sec:Pr}.
The Appendix contains some technical results.

\section{Adaptive procedure}\label{sec:Ad}

In this section we describe the adaptive procedure proposed in 
Galtchouk, Pergamenshchikov, 2006. We make use of the 
 standard trigonometric basis $(\phi_\zs{j})_\zs{j\ge 1}$ in $\cL_2[0,1]$, i.e.
\begin{equation}\label{sec:Ad.0}
\phi_1(x)=1\,,\quad
\phi_\zs{j}(x)=\sqrt{2}\,Tr_\zs{j}(2\pi [j/2]x)\,,\ j\ge 2\,,
\end{equation}
where the function $Tr_\zs{j}(x)=\cos(x)$ for even $j$ and
$Tr_\zs{j}(x)=\sin(x)$ for odd $j$; $[x]$ denotes the integer part of $x$.

To evaluate the error of estimation in the model \eqref{I.1} we will make use of the
empiric norm in the Hilbert space $\cL_\zs{2}[0,1]$,
generated by the design points $(x_\zs{j})_\zs{1\le j\le n}$ of model \eqref{I.1}. To this
end, for any functions $u$ and $v$ from $\cL_\zs{2}[0,1]$, we define the empiric
inner product
$$
(u\,,\,v)_\zs{n}=
\frac{1}{n}\sum^n_\zs{l=1}\,u(x_\zs{l})\,v(x_\zs{l})\,.
$$
Moreover, we will use this inner product for vectors in $\bbr^n$ as well, i.e.
if\\
$u=(u_\zs{1},\ldots,u_\zs{n})'$ and  $v=(v_\zs{1},\ldots,v_\zs{n})'$, then 
$$
(u\,,\,v)_\zs{n}=\frac{1}{n}u'v=
\frac{1}{n}\sum^n_\zs{l=1}\,u_\zs{l}\,v_\zs{l}\,.
$$
The prime denotes the transposition.

Notice that if $n$ is odd, then the functions $(\phi_\zs{j})_\zs{1\le j\le n}$ 
are orthonormal with respect to  this inner product, i.e. for any $1\le i,j\le n$,
\begin{equation}\label{sec:Ad.0-1}
(\phi_\zs{i}\,,\,\phi_\zs{j})_\zs{n}=
\frac{1}{n}\sum^n_\zs{l=1}\phi_\zs{i}(x_l)\phi_\zs{j}(x_l)={\bf Kr}_\zs{ij}\,,
\end{equation}
where ${\bf Kr}_\zs{ij}$ is Kronecker's symbol,  ${\bf Kr}_\zs{ij}=1$ if $i=j$ and 
 ${\bf Kr}_\zs{ij}=0$ for $i\ne j$.

\begin{remark}\label{Re.Ad.1}
Note that in the case of even $n$, the basis  \eqref{sec:Ad.0} is orthogonal and it is orthonormal except the 
$n$th function for which the normalizing constant should be changed.
The corresponding modifications of the formulas for even $n$  one can see in Galtchouk, Pergamenshchikov,2005.
To avoid these complications of formulas related to even $n$, we suppose $n$ to be odd.
\end{remark}

Thanks to this basis we pass to  the discrete Fourier transformation of model \eqref{I.1}:
\begin{equation}\label{sec:Ad.1}
\wh{\theta}_\zs{j,n}=\theta_\zs{j,n}+\frac{1}{\sqrt{n}}\xi_\zs{j,n}\,,
\end{equation}
where $\wh{\theta}_\zs{j,n}=(Y,\phi_\zs{j})_\zs{n}$, $Y=(y_\zs{1},\ldots,y_\zs{n})'$,
$\theta_\zs{j,n}=(S,\phi_\zs{j})_\zs{n}$
and
$$
\xi_\zs{j,n}=\frac{1}{\sqrt{n}}\sum^n_\zs{l=1}\sigma_l(S)\xi_l\phi_\zs{j}(x_l)\,.
$$

We estimate the function $S$ by the weighted least squares estimator
\begin{equation}\label{sec:Ad.2}
\wh{S}_\zs{\lambda}=\sum^n_\zs{j=1}\lambda(j)\wh{\theta}_\zs{j,n}\phi_\zs{j}\,,
\end{equation}
where the weight vector $\lambda=(\lambda(1),\ldots,\lambda(n))'$
belongs to some finite set $\Lambda$ from $[0,1]^n$ with $n\ge 3$.

Here we make use of the weight family $\Lambda$ introduced
in Galtchouk, Pergamenshchikov, 2008, i.e.
\begin{equation}\label{sec:Ad.3}
\Lambda\,=\,\{\lambda_\zs{\alpha}\,,\,\alpha\in\cA\}\,,
\quad
\cA=\{1,\ldots,k^*\}\times\{t_1,\ldots,t_m\}\,,
\end{equation}
where  $t_\zs{i}=i\ve$ and $m=[1/\ve^2]$. We suppose that the parameters 
$k^*\ge 1$ and $0<\ve\le 1$ are functions of $n$, i.e. $k^*= k^*_\zs{n}$ and 
$\ve=\ve_\zs{n}$, such that,
\begin{equation}\label{sec:Ad.3-1}
\left\{
\begin{array}{ll}
&\lim_\zs{n\to\infty}\,k^*_\zs{n}=+\infty\,,
\quad
\lim_\zs{n\to\infty}\,\dfrac{k^*_\zs{n}}{\ln n}=0\,,\\[4mm]
&\lim_\zs{n\to\infty}\,\ve_\zs{n}\,=\,0
\quad\mbox{and}\quad
\lim_\zs{n\to\infty}\,n^{\nu}\,\ve_\zs{n}\,=+\infty\,,
\end{array}
\right.
\end{equation}
for any $\nu>0$. For example, one can take for  $n\ge 3$
$$
\ve_\zs{n}=1/\ln n
\quad\mbox{and}\quad
k^*_\zs{n}=\ov{k}+\sqrt{\ln n}\,,
$$
where $\ov{k}$ is any nonnegative constant.

For each $\alpha=(\beta,t)\in\cA$ we define the weight vector
$\lambda_\zs{\alpha}=(\lambda_\zs{\alpha}(1),\ldots,\lambda_\zs{\alpha}(n))'$ as
\begin{equation}\label{sec:Ad.4}
\lambda_\zs{\alpha}(j)=\Chi_\zs{\{1\le j\le j_\zs{0}\}}+
\left(1-(j/\omega(\alpha))^\beta\right)\,
\Chi_\zs{\{ j_\zs{0}<j\le \omega(\alpha)\}}\,.
\end{equation}
Here $j_0=j_\zs{0}(\alpha)=\left[\omega(\alpha)\,\ve_\zs{n}\right]$ with
\begin{equation}\label{sec:Ad.4-1}
\omega(\alpha)=\ov{\omega}+(A_\zs{\beta}\,t)^{1/(2\beta+1)}n^{1/(2\beta+1)}\,,
 \end{equation}
where $\ov{\omega}$ is any nonnegative constant 
and
$$
A_\zs{\beta}=\frac{(\beta+1)(2\beta+1)}{\beta\pi^{2\beta}}\,.
$$

\begin{remark}\label{Re.Ad.2}
Note that the weighted least squares estimators 
 \eqref{sec:Ad.2} have been introduced by Pinsker, 1981, for continuous time optimal signal filtering in the gaussian noise. He proved that the mean-square asymptotic risk is minimized
by weighted least squares estimators with weights of type
  \eqref{sec:Ad.4}. Moreover he has found the sharp minimal value of the mean-square asymptotic
 risk, which was called later as the Pinsker constant.
 Nussbaum, 1985, used the same method with proper modification for efficient
estimation of the function $S$ of known smoothness
in the homogeneous gaussian model \eqref{I.1}, i.e. when
$\sigma_1(S)=\ldots=\sigma_\zs{n}(S)=1$ and $(\xi_\zs{j})_\zs{1\le j\le n}$ is i.i.d. $\cN(0,1)$ sequence. 
\end{remark}

 To choose weights from the set \eqref{sec:Ad.3} we minimize
the special cost function introduced by
 Galtchouk, Pergamenshchikov, 2007. This cost function is as follows
\begin{equation}\label{sec:Ad.5}
J_\zs{n}(\lambda)\,=\,\sum^n_\zs{j=1}\,\lambda^2(j)\wh{\theta}^2_\zs{j,n}\,-
2\,\sum^n_\zs{j=1}\,\lambda(j)\,\wt{\theta}_\zs{j,n}\,
+\,\rho \wh{P}_\zs{n}(\lambda)\,,
\end{equation}
where 
\begin{equation}\label{sec:Ad.6}
\wt{\theta}_\zs{j,n}=
\wh{\theta}^2_\zs{j,n}-\frac{1}{n}\wh{\varsigma}_\zs{n}
\quad\mbox{with}\quad
\wh{\varsigma}_\zs{n}=\sum^{n}_\zs{j=l_\zs{n}+1}
\wh{\theta}^2_\zs{j,n}
\end{equation}
and $l_\zs{n}=[n^{1/3}+1]$. The penalty term we define as
$$
\wh{P}_\zs{n}(\lambda)=\frac{|\lambda|^2 \wh{\varsigma}_\zs{n}}{n}\,,\quad
|\lambda|^2=\sum^n_\zs{j=1} \lambda^2(j)
\quad\mbox{and}\quad
\rho=\frac{1}{3+L_\zs{n}}\,,
$$
where $L_\zs{n}\ge 0$ is any slowly increasing sequence, i.e.
\begin{equation}\label{sec:Ad.6-1}
\lim_\zs{n\to\infty}\,L_\zs{n}=+\infty
\quad\mbox{and}\quad
\lim_\zs{n\to\infty}\,\frac{L_\zs{n}}{n^\nu}=0\,,
\end{equation}
for any $\nu>0$.

Finally, we set
\begin{equation}\label{sec:Ad.7}
\wh{\lambda}=\mbox{argmin}_\zs{\lambda\in\Lambda}\,J_\zs{n}(\lambda)
\quad\mbox{and}\quad
\wh{S}_\zs{*}=\wh{S}_\zs{\wh{\lambda}}\,.
\end{equation}

The goal of this paper is to study asymptotic (as $n\to\infty$) properties
of this estimation procedure.

\begin{remark}\label{Re.Ad.3}
Now we explain why does one choose the cost function in the form  \eqref{sec:Ad.5}.
Developing the empiric quadratic  risk for estimate \eqref{sec:Ad.2}, one obtains
$$
\|\wh{S}_\zs{\lambda}-S\|^2_\zs{n}=\sum^n_\zs{j=1}\,\lambda^2(j)\wh{\theta}^2_\zs{j,n}\,-
2\,\sum^n_\zs{j=1}\,\lambda(j)\,\wh{\theta}_\zs{j,n}\,\theta_\zs{j,n} +\|S\|^2_\zs{n}\,.
$$
It's natural to choose the weight vector  $\lambda$ for which this function reaches
the minimum. Since the last term on the right-hand part is independent of  $\lambda$, it can be droped and one has to minimize with respect to $\lambda$ the function equals the difference of two first terms on the right-hand part. It's clear that the minimization problem cann't be solved directly because the Fourier coefficients $(\theta_\zs{j,n}) $ are unknown.To overcome
this difficulty, we replace the product $\wh{\theta}_\zs{j,n}\,\theta_\zs{j,n}$
by its asymptotically unbiased estimator 
$\wt{\theta}_\zs{j,n}$ (see, Galtchouk, Pergamenshchikov, 2007, 2008). Moreover, to pay this
substitution, we introduce into the cost function the penalty term $\wh{P}_\zs{n}$ with a small coefficient $\rho>0$. The form of the penalty term is provided by the principal term of the quadratic risk for weighted least-squares estimator,
see Galtchouk, Pergamenshchikov, 2007, 2008.The coefficient $\rho>0$ means, that the penalty is small, because the estimator
 $\wt{\theta}_\zs{j,n}$ approximates in mean  the quantity
$\wh{\theta}_\zs{j,n}\,\theta_\zs{j,n}$ asymptotically, as $n\to\infty$.

Note that the principal difference between the procedure \eqref{sec:Ad.7} and the adaptive procedure proposed by Golubev, Nussbaum, 1993, for a homogeneous gaussian regression, consists in 
presence of the penalty term in the cost function \eqref{sec:Ad.5}. 
\end{remark}

\begin{remark}\label{Re.Ad.4}
As it was noted at Remark~\ref{Re.Ad.2}, Nussbaum, 1985, has shown that the weight coefficients
of type \eqref{sec:Ad.4} provide the asymptotic minimum of the mean-squared risk
at the regression function estimation problem for the homogeneous gaussian model  \eqref{I.1}, when
the smoothness of the function $S$ is known. In fact, to obtain an efficient estimator one needs to take a weighted least squares estimator \eqref{sec:Ad.2} with the weight vector
 $\lambda_\zs{\alpha}$, where the index $\alpha$ depends on smoothness of function $S$ and on
 coefficients $(\sigma_\zs{j}(S))_\zs{1\le j\le n}$,
(see \eqref{sec:Up.0-1} below), which are unknown in our case. For this reason, Galtchouk,
 Pergamenshchikov , 2007, 2008, have proposed to make use of the family of coefficients \eqref{sec:Ad.3}, which contains the weight vector providing the minimum of the mean-squared
risk. Moreover, they proposed the adaptive
 procedure  \eqref{sec:Ad.7} for which a non-asymptotic oracle inequality (see, Theorem~\ref{Th.M.1} below) was proved under some weak conditions on the  
coefficients $(\sigma_\zs{j}(S))_\zs{1\le j\le n}$. It is important to note that due the properties of the parametric family
\eqref{sec:Ad.3-1}, the secondary term in the oracle inequality is slowly increasing (slower than any degree of $n$).
\end{remark}

\section{Conditions}\label{sec:Co}

First we impose some conditions on unknown function $S$ in the model \eqref{I.1}.

Let   $\cC^{k}_\zs{per,1}(\bbr)$ be  the set of $1$-periodic
$k$ times differentiable $\bbr\to\bbr $ functions. We assume that
 $S$ belongs to the following set
\begin{equation}\label{sec:Co.1}
W^{k}_\zs{r}=\{f\in\cC^{k}_\zs{per,1}(\bbr)
\,:\,\sum_\zs{j=0}^k\,\|f^{(j)}\|^2\le r\}\,,
 \end{equation}
where $\|\cdot\|$ denotes the  norm in $\cL_\zs{2}[0,1]$, i.e.
\begin{equation}\label{sec:Co.2}
\|f\|^2=\int^1_\zs{0}f^2(t)\d t\,.
\end{equation}
Moreover, we suppose that $r>0$ and $k\ge 1$ are unknown parameters.

Note that, we can represent the set $W^{k}_\zs{r}$ as
an ellipse in $\cL_\zs{2}[0,1]$, i.e.
\begin{equation}\label{sec:Co.3}
W^{k}_\zs{r}=\{f\in\cL_\zs{2}[0,1]\,:\,
\sum_\zs{j=1}^\infty\,a_\zs{j}\theta^2_\zs{j}\le r\}\,,
 \end{equation}
where 
\begin{equation}\label{sec:Co.3-1}
\theta_\zs{j}=(f,\phi_\zs{j})=\int^1_\zs{0}f(t)\phi_\zs{j}(t)\d t
\end{equation}
and 
\begin{equation}\label{sec:Co.4}
a_\zs{j}=\sum^k_\zs{l=0}\|\phi^{(l)}_\zs{j}\|^2=
\sum^k_{i=0}(2\pi [j/2])^{2i}\,.
\end{equation}
Here $(\phi_\zs{j})_\zs{j\ge 1}$ is the trigonometric basis 
defined in \eqref{sec:Ad.0}.

Now we describe the conditions on the scale coefficients $(\sigma_\zs{j}(S))_\zs{j\ge 1}$.

\begin{itemize}
\item[$\H_\zs{1})$] {\em 
 $\sigma_\zs{j}(S)=g(x_\zs{j},S)$ for some unknown  function 
$g : [0,1]\times \cL_\zs{1}[0,1] \to \bbr_+$, which is
square integrable with respect to $x$ such that
\begin{equation}\label{sec:Co.5}
\lim_\zs{n\to\infty}\,\sup_\zs{S\in W^k_r}\,\left|
\,\frac{1}{n}\,\sum^n_\zs{j=1}\,g^2(x_\zs{j},S)\,-\,\varsigma(S)\,\right|\,=0\,,
\end{equation}
where $\varsigma(S):=\,\int_0^1\,g^2(x,S)\d x$. 
Moreover,
\begin{equation}\label{sec:Co.6}
g_\zs{*}=\inf_\zs{0\le x\le 1}\,
\inf_\zs{S\in W^k_r} g^2(x,S)\,>0
\end{equation}
and
\begin{equation}\label{sec:Co.6-1}
\sup_\zs{S\in W^k_r} \varsigma(S)<\infty\,.
\end{equation}
}
\item[$\H_\zs{2})$] {\em 
For any $x\in [0,1]$, the operator 
$g^2(x,\cdot)\,:\,\cC[0,1]\to \bbr$
 is differentiable in  
the Fr\'echet sense for any fixed function $f_\zs{0}$ from $\cC[0,1]$
, i.e.
 for any $f$ from some vicinity of $f_\zs{0}$ in $\cC[0,1]$,
$$
g^2(x,f)=g^2(x,f_\zs{0})+\L_\zs{x,f_\zs{0}}(f-f_\zs{0})+
\Upsilon(x,f_\zs{0},f)\,,
$$
where the Fr\'echet derivative
$\L_\zs{x,f_\zs{0}}\,:\,\cC[0,1]\to \bbr$
is a bounded linear operator
and 
the residual term $\Upsilon(x,f_\zs{0},f)$, for each $x\in [0,1]$, satisfies the following
property:
$$
\lim_\zs{\|f-f_\zs{0}\|_\zs{\infty}\to 0}
\frac{|\Upsilon(x,f_\zs{0},f)|}{\|f-f_\zs{0}\|_\zs{\infty}}=0\,,
$$
where $\|f\|_\zs{\infty}=\sup_\zs{0\le t\le 1} |f(t)|$.
}

\item[$\H_\zs{3})$] {\em
There exists some positive constant $C^*$ such that
for any function $S$ from $\cC[0,1]$ the operator 
$\L_\zs{x,S}$ defined in the condition $\H_\zs{2})$ 
satisfies the following inequality for any function $f$ from $\cC[0,1]$:
\begin{equation}\label{sec:Co.8}
|\L_\zs{x,S}(f)|
\le C^*
\left(
|S(x)f(x)|+|f|_\zs{1}+\|S\|\,\|f\|
\right)\,,
\end{equation}
where $|f|_\zs{1}=\int^1_\zs{0}|f(t)|\d t$.
}

\item[$\H_\zs{4})$] {\em The function 
$g_\zs{0}(\cdot)=g(\cdot,S_\zs{0})$ corresponding to $S_\zs{0}\equiv 0$
is  continuous on the interval $[0,1]$.
Moreover,
$$
\lim_\zs{\delta\to 0}\,
\sup_\zs{0\le x\le 1}\,
\sup_\zs{\|S\|_\zs{\infty}\le \delta}\,
|g(x,S)-g(x,S_\zs{0})|\,=\,0\,.
$$
}
\end{itemize}

\begin{remark}\label{Re.Co.1}
Let us explain the conditions $\H_\zs{1})$--$\H_\zs{4})$. In fact,this is the regularity conditions of the function $g(x,S)$ generating the scale coefficients $(\sigma_\zs{j}(S))_\zs{1\le j\le n}$.

Condition $\H_\zs{1})$ means that the function $g(\cdot,S)$ should be uniformly integrable
with respect to the first argument in the sens of convergence \eqref{sec:Co.5}.
Moreover, this function should be separated  from zero (see inequality \eqref{sec:Co.6})
and bounded on the class \eqref{sec:Co.1} (see inequality \eqref{sec:Co.6-1}). Boundedness away
from zero provides that the distribution of observations $(y_\zs{j})_\zs{1\le j\le n}$ isn't
degenerate in $\bbr^n$, and the boundedness means that the intensity of the noise vector should be finite, otherwise the estimation problem hasn't any sens.

 Conditions $\H_\zs{2})$ and $\H_\zs{3})$ mean that the function
$g(x,\cdot)$ is regular, at any fixed $0\le x\le 1$,  with respect to $S$ in the sens, that it is
differentiable in the Fr\`echet sens (see e.g., Kolmogorov, Fomin, 1989) and moreover the Fr\`echet derivative satisfies the growth condition given by the inequality \eqref{sec:Co.8}
which permits to consider the example \eqref{I.2}.

Last the condition $\H_\zs{4})$ is the usual uniform continuity the condition of the function $g(\cdot,\cdot)$ at the function $S_\zs{0}$.
\end{remark}

Now we give some examples of functions satisfying the conditions $\H_\zs{1})$-$\H_\zs{4})$.

We set
\begin{equation}\label{sec:Co.12}
g^2(x,S)=c_\zs{0}+c_\zs{1}x+c_\zs{2}S^2(x)+
c_\zs{3} \int^1_\zs{0}\,S^2(t) \d t
\end{equation}
with some coefficients $c_\zs{0}>0$,  $c_\zs{i}\ge 0$, $i=1,2,3$. 

In this case 
$$
\varsigma(S)=c_\zs{0}+\frac{c_\zs{1}}{2}+
(c_\zs{2}+c_\zs{3})
\int^1_\zs{0}\,S^2(t)\d t\,.
$$

The Fr\'echet derivative is given  by
$$
\L_\zs{x,S}(f)=2 S(x) f(x)+
2\int^1_\zs{0} S(t) f(t)\d t\,.
$$
It is easy to see that the function \eqref{sec:Co.12} satisfies the conditions $\H_\zs{1})$--$\H_\zs{4})$.
Moreover, the conditions $\H_\zs{1})$--$\H_\zs{4})$ are satisfied by any function of type
\begin{equation}\label{sec:Co.11}
g^2(x,S)=G(x,S(x))+\int^1_\zs{0}\,V(S(t))\d t\,,
\end{equation}
where the functions $G$ and $V$ satisfy the following the conditions:

\begin{itemize}
\item 
$G$ is a
 $[0,1]\times\bbr\to [c_\zs{0}\,,\,+\infty)$ function (with $c_\zs{0}>0$) such that
\begin{equation}\label{sec:Co.9}
\lim_\zs{\delta\to 0}\max_\zs{|u-v|\le \delta}
\sup_\zs{y\in\bbr} 
| G(u,y)-G(v,y)|=0
\end{equation}
and
\begin{equation}\label{sec:Co.10}
m_\zs{1}=
\sup_\zs{0\le x\le 1}
\sup_\zs{y\in\bbr}\,
\frac{|G_\zs{y}(x,y)|}{|y|}\,
<\infty\,;
\end{equation}

\item
 $V$ is a continuously differentiable  
$\bbr\to\bbr_\zs{+}$
function such that
$$
m_\zs{2}=
\,\sup_\zs{y\in\bbr}\,
\frac{|\dot{V}(y)|}{1+|y|}<\infty\,,
$$
where $\dot{V}(\cdot)$ is the derivative of $V$.

\end{itemize}

\noindent In this case 
$$
\varsigma(S)=\int^1_\zs{0}G(t,S(t))\d t+\int^1_\zs{0}\,V(S(t))\d t
$$
and
\begin{align*}
\left| n^{-1}\sum^n_\zs{j=1}\,g^2(x_\zs{j},S)-\varsigma(S)\right|&
\le
\sum^n_\zs{j=1}\,
\int^{x_\zs{j}}_\zs{x_\zs{j-1}}
\left|G(x_\zs{j},S(x_\zs{j}))-
G(t,S(t))
\right|\,\d t\\
&\le
\Delta_\zs{n}+
\sum^n_\zs{j=1}\,
\int^{x_\zs{j}}_\zs{x_\zs{j-1}}
\left|G(t,S(x_\zs{j}))-
G(t,S(t))
\right|\,\d t\,,
\end{align*}
where $\Delta_\zs{n}=\max_\zs{|u-v|\le 1/n}\sup_\zs{y\in\bbr}|G(u,y)-G(v,y)|$.
Now to estimate the last term in this inequality note that
$$
G(t,S(x_\zs{j}))-
G(t,S(t))=\int^{x_\zs{j}}_\zs{t}\,G_\zs{y}(t,S(z))\,\dot{S}(z)\,\d z\,.
$$
Therefore, from the condition \eqref{sec:Co.10} we get 
$$
|G(t,S(x_\zs{j}))-
G(t,S(t))|\le m_\zs{1}\int^{x_\zs{j}}_\zs{x_\zs{j-1}}\,|S(z)|\,|\dot{S}(z)|\,\d z\,,
$$
and through the Bounyakovskii-Cauchy-Schwarz inequality,  for any $S\in W^k_\zs{r}$,
\begin{align*}
\left| n^{-1}\sum^n_\zs{j=1}\,g^2(x_\zs{j},S)-\varsigma(S)\right|&
\le
\sum^n_\zs{j=1}\,
\int^{x_\zs{j}}_\zs{x_\zs{j-1}}
\left|G(x_\zs{j},S(x_\zs{j}))-
G(t,S(t))
\right|\,\d t\\
&\le
\Delta_\zs{n}+
\frac{m_\zs{1}}{n}\,\int^1_\zs{0}\,|S(t)|\,|\dot{S}(t)|\d t\\
&\le
\Delta_\zs{n}+
\frac{m_\zs{1}}{n}\,
\|S\|\,\|\dot{S}\|\,
\le \Delta_\zs{n}+\frac{m_\zs{1}}{n}\, r\,.
\end{align*}
Now, the condition \eqref{sec:Co.9} implies $\H_\zs{1})$.

Moreover, the Fr\'echet derivative in this case is given  by
$$
\L_\zs{x,S}(f)=G_\zs{y}(x,S(x))f(x)+
\int^1_\zs{0}\dot{V}(S(t)) f(t)\d t\,.
$$
One can check directly that this operator satisfies the inequality
  \eqref{sec:Co.8} with $C^*=m_\zs{1}+m_\zs{2}$.


\medskip
\vspace{5mm}
\section{ Main results}\label{sec:Ma}

Denote by $\cP_\zs{n}$  the family of distributions $p$ in $\bbr^n$ of the vectors
$(\xi_\zs{1},\ldots,\xi_\zs{n})'$
in the model \eqref{I.1} such that the components $\xi_\zs{j}$ are jointly independent,
centered with unit variance and
\begin{equation}\label{sec:Ma.0}
\max_\zs{1\le k\le n}\,
\E\, \xi_\zs{k}^4\le l^*_\zs{n}\,,
\end{equation}
where $l^*_\zs{n}\ge 3$ is slowly increasing sequence, that is
it satisfies the property \eqref{sec:Ad.6-1}.

It is easy to see that, for any $n\ge 1$, the centered gaussian distribution in $\bbr^n$ with
unit covariation matrix belongs to the family $\cP_\zs{n}$. We will denote by $q$ this gaussian distribution.



For any estimator $\wh{S}$ we define the following quadratic risk
\begin{equation}\label{sec:Ma.1}
\cR_\zs{n}(\wh{S},S)=\sup_\zs{p\in\cP_\zs{n}}\E_\zs{S,p}\|\wh{S}-S\|^2_\zs{n}\,,
\end{equation}
where $\E_\zs{S,p}$ is the expectation with respect to the distribution $\P_\zs{S,p}$
of the observations $(y_\zs{1},\ldots,y_\zs{n})$ with the fixed function $S$ and
the fixed  distribution $p\in\cP_\zs{n}$ of random variables $(\xi_\zs{j})_\zs{1\le j\le n}$ in the model \eqref{I.1}.

Moreover, to make the risk independent of the design points, 
 in this paper we will make use of the risk with respect to the usual norm in $\cL_\zs{2}[0,1]$ \eqref{sec:Co.2} also,
i.e.
\begin{equation}\label{sec:Ma.1-1}
\cT_\zs{n}(\wh{S},S)=\sup_\zs{p\in\cP_\zs{n}}\E_\zs{S,p}\|\wh{S}-S\|^2\,.
\end{equation}

If an estimator $\wh{S}$ is defined only at the design points $(x_\zs{j})_\zs{1\le j\le n}$, then we extend it as step function onto the interval $[0,1]$ by setting $\wh{S}(x)=T(\wh{S}(x))$,
for all $0\le x\le 1$, where
\begin{equation}\label{sec:Ma.1-2}
T(f)(x)=f(x_\zs{1})\Chi_\zs{[0,x_{1}]}(x)
+\sum_{k=2}^n\,f(x_k)\Chi_\zs{(x_{k-1},x_k]}(x)\,.
\end{equation}

In Galtchouk, Pergamenshchikov, 2007, 2008
the following non-asymptotic oracle inequality has been shown for the procedure 
\eqref{sec:Ad.7} .
\begin{theorem}\label{Th.M.1}
Assume that in the model \eqref{I.1} the function
$S$ belongs to  $W_\zs{r}^{1}$.
Then, for any odd $n\ge 3$ and $r>0$,
 the estimator $\wh{S}_\zs{*}$ satisfies the following 
oracle inequality
\begin{equation}\label{sec:Ma.2}
\cR_\zs{n}(\wh{S}_\zs{*},S)\,
\le\,
\frac{1+3\rho-2\rho^2}{1-3\rho}
\,\min_\zs{\lambda\in\Lambda}\,
\cR_\zs{n}(\wh{S}_\zs{\lambda},S)\,
+\,\frac{1}{n}\,\cB_\zs{n}(\rho)\,,
\end{equation}
where the function $\cB_\zs{n}(\rho)$ is such that, for any $\nu>0$,
\begin{equation}\label{sec:Ma.3}
\lim_\zs{n\to\infty}\,
\frac{\cB_\zs{n}(\rho)}{n^\nu}\,=\,0\,.
\end{equation}
\end{theorem}

\begin{remark}\label{Re.M.1}
Note that in Galtchouk, Pergamenshchikov, 2007, 2008, the oracle inequality is proved
for the model \eqref{I.1}, where the random variables $(\xi_\zs{j})_\zs{1\le j\le n}$
are independent identically distributed. In fact, the result and the proof are true for 
independent random variables which are not identically distributed, i.e. for any
distribution of the random vector $(\xi_\zs{1},\ldots,\xi_\zs{n})'$  from $\cP_\zs{n}$.
\end{remark}

Now we formulate the main asymptotic results. To this end, 
for any function  $S\in W^k_\zs{r}$, we set
\begin{equation}\label{sec:Ma.4}
\gamma_k(S)\,=\,\Gamma^*_k\,r^{\frac{1}{2k+1}}\,(\varsigma(S))^{2k/(2k+1)}\,,
\end{equation}
where
$$
\Gamma^*_k=(2k+1)^{\frac{1}{2k+1}}\left(k/(\pi\,(k+1))\right)^{2k/(2k+1)}\,.
$$

It is well known (see e.g., Nussbaum, 1985)
that the optimal rate of convergence is $n^{2k/(2k+1)}$ when 
the risk is taken uniformly over $W^k_\zs{r}$.

\begin{theorem}\label{Th.M.2}
Assume that in the model \eqref{I.1} the sequence
$(\sigma_\zs{j}(S))$ fulfills the condition 
$\H_1)$. Then the estimator $\wh{S}_\zs{*}$ from \eqref{sec:Ad.7}
 satisfies the inequalities
\begin{equation}\label{sec:Ma.5}
\limsup_{n\to\infty}\,
n^{\frac{2k}{2k+1}}\,
\sup_\zs{S\in W^k_\zs{r}}\,
\frac{\cR_\zs{n}(\wh{S}_\zs{*},S)}{\gamma_k(S)}\,\le 1
\end{equation}
and
\begin{equation}\label{sec:Ma.5-1}
\limsup_{n\to\infty}\,
n^{\frac{2k}{2k+1}}\,
\sup_\zs{S\in W^k_\zs{r}}\,
\frac{\cT_\zs{n}(\wh{S}_\zs{*},S)}{\gamma_k(S)}\,\le 1\,.
\end{equation}
\end{theorem}
The following result  gives the sharp lower bound for risk \eqref{sec:Ma.1} and show that
$\gamma_k(S)$ is the \mbox{\it Pinsker} \mbox{\rm constant}.
\begin{theorem}\label{Th.M.3}
Assume that in the model \eqref{I.1}  the sequence
$(\sigma_\zs{j}(S))$
satisfies the conditions $\H_2)$-- $\H_4)$. Then the risks
\eqref{sec:Ma.1} and \eqref{sec:Ma.1-1}
admit the following asymptotic lower bounds
\begin{equation}\label{sec:Ma.6}
\liminf_{n\to\infty}\,n^{\frac{2k}{2k+1}}\,\inf_{\wh{S}_\zs{n}}\,\sup_\zs{S\in W^k_\zs{r}}\,
\frac{\cR_\zs{n}(\wh{S}_\zs{n},S)}{\gamma_k(S)}\,
\ge 1
\end{equation}
and
\begin{equation}\label{sec:Ma.6-1}
\liminf_{n\to\infty}\,n^{\frac{2k}{2k+1}}\,\inf_{\wh{S}_\zs{n}}\,\sup_\zs{S\in W^k_\zs{r}}\,
\frac{\cT_\zs{n}(\wh{S}_\zs{n},S)}{\gamma_k(S)}\,
\ge 1\,.
\end{equation}
\end{theorem}

\begin{remark}\label{Re.M.2}
Note that in Galtchouk, Pergamenshchikov, 2005,
an asymptotically efficient estimator has been constructed and results
similar to Theorems \ref{Th.M.2} and \ref{Th.M.3} were claimed for the model  \eqref{I.1}.
In fact the upper bound is true there under some additional condition on the smoothness of
the function $S$, i.e. on the parameter $k$. In the cited paper this additional condition
is not formulated since erroneous inequality $(A.6)$. To avoid using  this inequality
we modify the estimating procedure by introducing the penalty term 
$\rho\, \wh{P}_\zs{n}(\lambda)$ in the cost function \eqref{sec:Ad.5}. By this way we remove
all additional conditions on the smoothness parameter $k$.
\end{remark}

\begin{remark}\label{Re.M.3}
In fact to obtain the non-asymptotic oracle inequality \eqref{sec:Ma.2}, it isn't necessary to make use of equidistant design points and the trigonometric basis. One may take any design points (deterministic or random) and any orthonormal basis satisfying \eqref{sec:Ad.0-1}. But
to obtain the property \eqref{sec:Ma.3} one needs to impose some technical conditions (see
Galtchouk, Pergamenshchikov, 2008).

Note that the results of Theorem~\ref{Th.M.2} and Theorem~\ref{Th.M.3} are based on equidistant design points and the trigonometric basis.


\end{remark}

\medskip
\section{ Upper bound}\label{sec:Up}

In this section we prove Theorem~\ref{Th.M.2}.  To this end we will make use of the oracle inequality
\eqref{sec:Ma.2}. We have to find an estimator from the family
\eqref{sec:Ad.2}-\eqref{sec:Ad.3} for which we can show the upper bound \eqref{sec:Ma.5}. We start with the construction
of  such an estimator. First we put
\begin{equation}\label{sec:Up.0}
\wt{l}_\zs{n}=\inf\{i\ge 1\,:\,i\ve\ge \ov{r}(S)\}\wedge m
\quad\mbox{and}\quad
\ov{r}(S)=r/\varsigma(S)\,,
\end{equation}
where $a\wedge b=\min(a,b)$.

Then we choose an index from the set $\cA_\ve$ as
\begin{equation}\label{sec:Up.0-1}
\wt{\alpha}=(k,\wt{t}_\zs{n})\,,
\end{equation}
where $k$ is the parameter of the set $W^k_\zs{r}$ and $\wt{t}_\zs{n}=\wt{l}_\zs{n}\ve$.
Finally, we set
\begin{equation}\label{sec:Up.0-2}
\wt{S}=\wh{S}_\zs{\wt{\lambda}}
\quad\mbox{and}\quad
\wt{\lambda}=\lambda_\zs{\wt{\alpha}}\,.
\end{equation}
Now we show the upper bound \eqref{sec:Ma.5} for this estimator.
\begin{theorem}\label{Th.U.1} Assume that the condition $\H_1)$ holds. Then 
\begin{equation}\label{sec:Up.1}
\limsup_{n\to\infty}\,n^{\frac{2k}{2k+1}}\,
\sup_{S\in W^k_r}\,
\frac{\cR_\zs{n}(\wt{S},S)}{\gamma_k(S)}\,\le 1\,.
\end{equation}
\end{theorem}
\begin{remark}\label{Re.U.1}
Note that the estimator
$\wt{S}$ belongs to the family \eqref{sec:Ad.2}-\eqref{sec:Ad.3}, but we can't use directly
 this estimator because the parameters $k$, $r$ and $\ov{r}(S)$ are unknown. We can use
this upper bound only through the oracle inequality \eqref{sec:Ma.2} proved for 
procedure \eqref{sec:Ad.7}.
\end{remark}

Now Theorem~\ref{Th.M.1} and Theorem~\ref{Th.U.1} imply the upper bound
\eqref{sec:Ma.5}. To obtain the upper bound \eqref{sec:Ma.5-1} we need the following auxiliary result.

\begin{lemma}\label{sec:Le.U.1}
For any $0<\delta<1$ and  any estimate $\wh{S}_\zs{n}$ of $S\in W^k_\zs{r}$,
$$
\|\wh{S}_\zs{n}-S\|_\zs{n}^2\ge(1-\delta)\|T_\zs{n}(\wh{S})-S\|^2\,
-\,(\delta^{-1}\,-\,1)\,r/n^2\,,
$$
where the function
$T_\zs{n}(\wh{S})(\cdot)$ is defined in \eqref{sec:Ma.1-2}.
\end{lemma}
\noindent Proof of this Lemma is given in Appendix~\ref{subsec:A.1}.

Now inequality \eqref{sec:Ma.5} and this lemma imply the upper bound \eqref{sec:Ma.5-1}.
Hence Theorem~\ref{Th.M.2}.

\medskip
\section{ Lower bound}\label{sec:Lo}

In this section we give the main steps of proving the lower bounds \eqref{sec:Ma.6} and
\eqref{sec:Ma.6-1}. In common, we follow the same scheme as Nussbaum, 1985.
We begin with minorizing the minimax risk by a bayesian one constructed on a parametric functional family introduced in Section~\ref{subsec:Fa} ( see \eqref{Fa.4}) and using the
prior distribution \eqref{Fa.5}.
Further, a special modification of the van Trees inequality (see, Theorem~\ref{Th.Tr.1}) yields
a lower bound for the bayesian risk depending on the chosen prior distribution, of course.
Finally, in section~\ref{subsec:Br}, we choose parameters of the prior distribution
(see \eqref{Fa.5}) providing the maximal value of the lower bound for the bayesian risk. This
value coincides with the Pinsker constant as it is shown in Section~\ref{subsec:Th.M.3}.

\subsection{ Lower bound for parametric heteroscedastic regression models} \label{subsec:Tr}

Let $(\bbr^n,\cB(\bbr^n),\P_\zs{\vartheta},\vartheta\in\Theta\subseteq\bbr^l)$ 
be a 
 statistical model relative to the observations $(y_\zs{j})_\zs{1\le j\le n}$ 
governed by the regression equation
\begin{equation}\label{subsec:Tr.1}
y_\zs{j}\,=\,S_\zs{\vartheta}(x_\zs{j})\,+\,\sigma_\zs{j}(\vartheta)\,\xi_\zs{j}\,,
\end{equation}
where $\xi_1,\ldots,\xi_\zs{n}$ are i.i.d. $\cN(0,1)$ random variables,
 $\vartheta=(\vartheta_1,\ldots,\vartheta_l)^\prime$ is a unknown parameter vector,
 $S_\zs{\vartheta}(x)$ is a unknown (or known) function 
and $\sigma_\zs{j}(\vartheta)=g(x_\zs{j},S_\vartheta)$, with the function
$g(x,S)$ defined in the condition $\H_\zs{1})$. Assume that a prior distribution
$\mu_\zs{\vartheta}$ of the parameter $\vartheta$ in $\bbr^l$ is defined by the density
 $\Phi(\cdot)$ of
the following form
$$
\Phi(z)\,=\,
\Phi(z_1,\ldots,z_l)=\prod_{i=1}^l\varphi_\zs{i}(z_\zs{i})\,,
$$
where $\varphi_\zs{i}$ is a continuously differentiable bounded density on $\bbr$ with
$$
I_\zs{i}=\int_{\bbr}\frac{\dot{\varphi}^2_\zs{i}(u)}{\varphi_\zs{i}(u)}\d u
<\infty\,.
$$
Let $\tau(\cdot)$ be a continuously differentiable $\bbr^l\to \bbr$ function
such that, for any $1\le i\le l$,
\begin{equation}\label{subsec:Tr.1-1}
\lim_\zs{|z_\zs{i}|\to\infty}\,\tau(z)\,\varphi_\zs{i}(z_\zs{i})=0
\quad\mbox{and}\quad
\int_{\bbr^l}\,\left|
\tau^{\prime}_\zs{i}(z)
\right|\,
\Phi(z)\d z<\infty\,,
\end{equation}
where 
$$ 
\tau^{\prime}_\zs{i}(z)=
(\partial/\partial z_\zs{i})\,\tau(z)\,.
$$
Let $\wh{\tau}_\zs{n}$ be an estimator of $\tau(\vartheta)$ based on
 observations $(y_\zs{j})_\zs{1\le j\le n}$.
 For any $\cB(\bbr^n\times\bbr^l)$ -
measurable integrable function 
$G(x,z), x\in\bbr^n, z\in\bbr^l$, 
we set
$$
\wt{\E}\,G(Y,\vartheta)\,=\,\int_{\bbr^l}\,\E_\zs{z}\,G(Y,z)\,
\Phi(z)\,\d z\,,
$$
where $\E_\zs{\vartheta}$ is the expectation with respect to the distribution 
$\P_\zs{\vartheta}$ of the vector $Y=(y_\zs{1},\ldots,y_\zs{n})$.
 Note that in this case
$$
\E_\zs{\vartheta}\,G(Y,\vartheta)\,=\,\int_\zs{\bbr^n}\,G(v,\vartheta)\,
f(v,\vartheta)\,\d v\,,
$$
where
\begin{equation}\label{subsec:Tr.2}
f(v,z)=\prod^n_\zs{j=1}\,\frac{1}{\sqrt{2\pi} \sigma_\zs{j}(z)}
\exp\left\{-\,\frac{(v_\zs{j}-S_\zs{z}(x_\zs{j}))^2}
{2\sigma_\zs{j}^2(z)}\right\}\,.
\end{equation}

We prove the following result.
\begin{theorem}\label{Th.Tr.1}
Assume that the conditions $\H_1)-\H_2)$ hold.
Moreover, assume that 
 the function $S_\zs{z}(\cdot)$ with $z=(z_\zs{1},\ldots,z_\zs{l})'$ is
uniformly over $0\le x\le 1$
 differentiable with respect to
$z_\zs{i},\ 1\le i\le l$, i.e. 
for any $1\le i\le l$
there exists
a function $S^{\prime}_\zs{z,i}\in \cC[0,1]$
such that
\begin{equation}\label{subsec:Tr.2-1}
\lim_\zs{h\to 0}\,
\max_\zs{0\le x\le 1}
\left|\left(
S_\zs{z+h\e_\zs{i}}(x)-S_\zs{z}(x)-S^{\prime}_\zs{z,i}(x)h\right)/h
\right|\,=0\,,
\end{equation}
where $\e_\zs{i}=(0,....,1,...,0)'$,  all coordinates are $0$, except the i-th equals to $1$ .
Then for any square integrable estimator $\wh{\tau}_\zs{n}$ of 
$\tau(\vartheta)$ and any $1\le i\le l$,
\begin{equation}\label{subsec:Tr.3}
\wt{\E}\left(\wh{\tau}_\zs{n}-\tau(\vartheta)\right)^2\ge
\frac{\ov{\tau}_\zs{i}^2}{F_\zs{i}\,+\,B_\zs{i}\,+\,I_\zs{i}}\,,
\end{equation}
where $\ov{\tau}_\zs{i}=
\int_{\bbr^l}\,
\tau^{\prime}_\zs{i}(z)\,
\Phi(z)\d z$, $F_\zs{i}=\sum_{j=1}^n
\int_\zs{\bbr^l}
\,(S^{\prime}_\zs{z,i}(x_\zs{j})/\sigma_\zs{j}(z))^2\,
\Phi(z)\d z$
and
$$
 B_\zs{i}=\,\frac{1}{2}\,\sum^n_\zs{j=1}\int_\zs{\bbr^l}
\frac{\wt{\L}^2_\zs{i}(x_\zs{j},S_\zs{z})}{\sigma^4_\zs{j}(S_\zs{z})}
\Phi(z)\d z\,, 
$$
$\wt{\L}_\zs{i}(x,z)=\L_\zs{x,S_\zs{z}}(S^{\prime}_\zs{z,i})$,
the operator $\L_\zs{x_,S}$ is defined in 
the condition $\H_2)$.
\end{theorem}
\noindent {\bf  Proof} is given in Appendix~\ref{subsec:A.2}.
\begin{remark}\label{Re.Tr.1}
Note that the inequality \eqref{subsec:Tr.3} is  some modification of the van Trees
inequality (see, Gill, Levit, 1995) adapted to the model \eqref{subsec:Tr.1}.
\end{remark}
\medskip

\subsection{ Parametric family of kernel functions}\label{subsec:Fa}

In this section we define and study some special parametric family of kernel function  
which will be used to prove the sharp lower bound \eqref{sec:Ma.6}.

Let us begin by kernel functions. We fix $\eta >0$ and we set 
\begin{equation}\label{Fa.1}
\chi_\zs{\eta}(x)
=\eta^{-1}\int_{\bbr}\Chi_\zs{(|u|\le 1-\eta)}\,V\left(\frac{u-x}{\eta}\right)\,\d u\,,
\end{equation}
where $\Chi_A$ is the indicator of a set $A$, the kernel $V\in\cC^{\infty}(\bbr)$  is 
such that
$$
V(u)=0
\quad\mbox{for}\quad |u|\ge 1
\quad\mbox{and}\quad 
\int^1_{-1}\,V(u)\,\d u=1\,.
$$
It is easy to see that the function $\chi_\zs{\eta}(x)$ possesses the properties :
\begin{align*}
&0\le \chi_\zs{\eta}\le 1\,,\quad \chi_\zs{\eta}(x)=1
\quad\mbox{for}\quad |x|\le 1-2\eta \quad\mbox{and}\\ 
&\chi_\zs{\eta}(x)=0
\quad\mbox{for}\quad |x|\ge 1\,.
\end{align*}
 Moreover,
for any $c>0$ and $\nu\ge 0$
\begin{equation}\label{Fa.2}
\lim_\zs{\eta\to 0}\,\sup_\zs{f\,:\,\|f\|_\zs{\infty}\le c}\,
\left|
\int_{\bbr}\,f(x)\,\chi_\zs{\eta}^\nu(x)\d x-\int_{-1}^{1}f(x) \d x
\right|
=0\,.
\end{equation}

We divide the interval $[0,1]$ into $M$ equal subintervals of length $2h$ and on each of them
we construct a  kernel-type function 
which equals to zero
at the boundary of the subinterval together with all derivatives. 

It provides that the Fourier partial sums with respect
to the trigonometric basis in $\cL_\zs{2}[-1,1]$ give a natural parametric approximation to the
function on each subinterval. 

Let  $(e_\zs{j})_\zs{j\ge 1}$ be the trigonometric basis in $\cL_2[-1,1]$, i.e.
\begin{equation}\label{Fa.3}
e_\zs{1}=1/\sqrt{2}\,,\quad
e_\zs{j}(x)=\,Tr_\zs{j}\left(\pi [j/2] x\right)\,,\ j\ge 2\,,
\end{equation}
where the functions $(Tr_\zs{j})_\zs{j\ge 2}$ are defined in \eqref{sec:Ad.0}.

Now,
 for any array $z=\{(z_\zs{m,j})_\zs{1\le m\le M_\zs{n}\,,\,1\le j\le N_\zs{n}}\}$  we define
the following function
\begin{equation}\label{Fa.4}
S_\zs{z,n}(x)=\sum_{m=1}^{M_\zs{n}}\sum_{j=1}^{N_\zs{n}}\,z_\zs{m,j}\,D_\zs{m,j}(x)\,,
\end{equation}
where $D_\zs{m,j}(x)=e_\zs{j}\left(v_m(x)\right)\chi_\zs{\eta}\left(v_m(x)\right)$,
$$
v_m(x)=(x-\wt{x}_m)/h_\zs{n}\,,
\quad\wt{x}_m= 2mh_\zs{n}
\quad\mbox{and}\quad
M_\zs{n}=\left[1/(2h_\zs{n})\right]-1\,.
$$

We assume that the sequences 
$(N_\zs{n})_\zs{n\ge 1}$ and $(h_\zs{n})_\zs{n\ge 1}$,
satisfy the following conditions.

$\A_\zs{1})$
{\em The sequence $N_\zs{n}\to\infty$ as $n\to\infty$ and for any $\nu>0$
$$
\lim_\zs{n\to\infty}\,N^{\nu}_\zs{n}/n\,=\,0\,.
$$
Moreover, there exist $0<\delta_\zs{1}<1$ and $\delta_\zs{2}>0$
such that 
$$
h_\zs{n}=\,\O(n^{-\delta_\zs{1}})
\quad\mbox{and}\quad
h^{-1}_\zs{n}=\,\O(n^{\delta_\zs{2}})
\quad\mbox{as}\quad
n\to \infty\,.
$$
}
To define a prior distribution on the family of arrays,
we choose the following random array 
$\vartheta=\{(\vartheta_\zs{m,j})_\zs{ 1\le m\le M_\zs{n}\,,\, 1\le j\le N_\zs{n}}\}$ 
with
\begin{equation}\label{Fa.5}
\vartheta_\zs{m,j}\,=\,t_\zs{m,j}\,\zeta_\zs{m,j}\,,
\end{equation}
where  $(\zeta_\zs{m,j})$ are i.i.d. $\cN(0,1)$ random variables and 
$(t_\zs{m,j})$
are some nonrandom positive coefficients. We make use of gaussian variables since they
possess the minimal Fisher information and therefore maximize the lower bound
\eqref{subsec:Tr.3}.
 We set
\begin{equation}\label{Fa.6}
t^*_\zs{n}=\,\max_\zs{1\le m\le M_\zs{n}}
\sum^{N_\zs{n}}_\zs{j=1}\,t_\zs{m,j}\,.
\end{equation}
We assume that the coefficients $(t_\zs{m,j})_\zs{ 1\le m\le M_\zs{n}\,,\, 1\le j\le N_\zs{n}}$
satisfy the following conditions.

$\A_\zs{2})$
{\em There exists  a sequence of positive numbers $(d_\zs{n})_\zs{n\ge 1}$ such that
\begin{equation}\label{Fa.7}
\lim_\zs{n\to\infty}
\frac{d_\zs{n}}{h_\zs{n}^{2k-1}}\,
\sum^{M_\zs{n}}_\zs{m=1}\sum^{N_\zs{n}}_\zs{j=1}\,t^2_\zs{m,j}\,
j^{2(k-1)}=0\,,
\quad
\lim_\zs{n\to\infty}\,\sqrt{d_\zs{n}}\,t^*_\zs{n}=0\,,
\end{equation}
moreover, for any $\nu> 0$,
$$
\lim_\zs{n\to\infty}\,n^{\nu}\,\exp\{{-d_\zs{n}/2}\}=0\,.
$$
}

$\A_\zs{3})$
{\em For some $0<\epsilon<1$
$$
\limsup_\zs{n\to\infty}
\frac{1}{h_\zs{n}^{2k-1}}\,\sum^{M_\zs{n}}_\zs{m=1}\sum^{N_\zs{n}}_\zs{j=1}\,t^2_\zs{m,j}\,
j^{2k}\,
\le (1-\epsilon)r
\left(\frac{2}{\pi}\right)^{2k}
\,.
$$
}

$\A_\zs{4})$
{\em There exists $\epsilon_\zs{0}>0$ such that
$$
\lim_\zs{n\to\infty}
\frac{1}{h_\zs{n}^{4k-2+\epsilon_\zs{0}}}\,
\sum^{M_\zs{n}}_\zs{m=1}\sum^{N_\zs{n}}_\zs{j=1}\,t^4_\zs{m,j}\,j^{4k}\,=0\,.
$$
}

\begin{prop}\label{sec:Pr.Fa.1}
Let the conditions $\A_\zs{1})$--$\A_\zs{2})$. Then,
for any $\nu>0$ and for any $\delta>0$,
$$
\lim_\zs{n\to\infty}\,n^\nu\,\max_\zs{0\le l\le k-1}
\P\left(\|S^{(l)}_\zs{\vartheta,n}\|>\delta\right)=0\,.
$$
\end{prop}

\begin{prop}\label{sec:Pr.Fa.2}
Let the conditions
$\A_\zs{1})$--$\A_\zs{4})$. Then, for any $\nu>0$,
$$
\lim_\zs{n\to\infty}\,n^\nu\,
\P(S_\zs{\vartheta,n}\notin W^{k}_\zs{r})\,=0\,.
$$
\end{prop}

\medskip

\begin{prop}\label{sec:Pr.Fa.3}
Let the conditions
$\A_\zs{1})$--$\A_\zs{4})$. Then, for any $\nu>0$,
$$
\lim_\zs{n\to\infty}\,n^\nu\,
\E\,\|S_\zs{\vartheta,n}\|^2\,
\left(
\Chi_\zs{\{S_\zs{\vartheta,n}\notin W^{k}_\zs{r}\}}\,+\,
\Chi_\zs{\Xi^c_\zs{n}}
\right)
=0\,.
$$
\end{prop}
\medskip

\begin{prop}\label{sec:Pr.Fa.4}
Let the conditions
$\A_\zs{1})$--$\A_\zs{4})$. Then for any function $g$
satisfying the conditions \eqref{sec:Co.6} and $\H_\zs{4})$ 
$$
\lim_\zs{n\to\infty}\,\sup_\zs{0\le x\le 1}\,
\E\,\left|\,g^{-2}(x,S_\zs{\vartheta,n})-g^{-2}_\zs{0}(x)\right|=0\,.
$$
\end{prop}

Proofs of Propositions~\ref{sec:Pr.Fa.1}--\ref{sec:Pr.Fa.4}
are given in Appendix.
\medskip

\subsection{Bayes risk}\label{subsec:Br}

Now we will obtain the lower bound for the bayesian risk that yields the lower bound
\eqref{sec:Ma.6-1} for the minimax risk.

We make use of the sequence of random functions
$(S_\zs{\vartheta,n})_\zs{n\ge 1}$ defined in \eqref{Fa.4}-\eqref{Fa.5}
with the coefficients $(t_\zs{m,j})$ satisfying the conditions $\A_\zs{1})$--$\A_\zs{4})$
which will be chosen later.

For any estimator $\wh{S}_\zs{n}$ we introduce now the corresponding Bayes risk
\begin{equation}\label{sec:Lo.3}
\cE_\zs{n}(\wh{S}_\zs{n})=
\int_\zs{\bbr^l}\,\E_\zs{S_\zs{z,n},q}\|\wh{S}_\zs{n}-S_\zs{z,n}\|^2\,\mu_\zs{\vartheta}(\d z)\,,
\end{equation}
where the kernel family $(S_\zs{z,n})$ is defined in \eqref{Fa.4}, $\mu_\zs{\vartheta}$ denotes 
 the distribution of the random array
 $\vartheta$ defined by \eqref{Fa.5} in $\bbr^l$ with $l= M_\zs{n}N_\zs{n}$.

We remember that $q$ is a centered gaussian distribution in $\bbr^n$ with unit covariation
 matrix.

First of all, we replace the functions
$\wh{S}_\zs{n}$ and $S$ by their Fourier series with respect to the basis
$$
\wt{e}_\zs{m,i}(x)=(1/\sqrt{h})\,e_\zs{i}\left(v_m(x)\right)\,
\Chi_\zs{\left(|v_m(x)|\le 1\right)}\,.
$$
By making use of this basis
we can estimate the norm $\|\wh{S}_\zs{n}-S_\zs{z,n}\|^2$ from below as
$$
\|\wh{S}_\zs{n}-S_\zs{z,n}\|^2 \ge
 \sum_{m=1}^{M_\zs{n}}\sum_{j=1}^{N_\zs{n}}\,(\wh{\tau}_\zs{m,j}\,-\,
\tau_\zs{m,j}(z))^2\,,
$$
where
$$
\wh{\tau}_\zs{m,j}=\int_0^1\,\wh{S}_\zs{n}(x)\wt{e}_\zs{m,j}(x)\d x
\quad\mbox{and}\quad 
\tau_\zs{m,j}(z)=\int_0^1\,S_\zs{z,n}(x)\wt{e}_\zs{m,j}(x)\,\d x\,.
$$
Moreover, from the definition \eqref{Fa.4}  one gets
$$
\tau_\zs{m,j}(z)
=\sqrt{h}\sum_{i=1}^{N_\zs{n}}\,z_\zs{m,i}\int_{-1}^1\,e_\zs{i}(u)e_\zs{j}(u)\chi_\zs{\eta}(u)\,\d u\,.
$$
It is easy to see that the functions $\tau_\zs{m,j}(\cdot)$
satisfy the condition \eqref{subsec:Tr.1-1} for gaussian prior densities. In this case
(see the definition in \eqref{subsec:Tr.3}) we have
$$
\ov{\tau}_\zs{m,j}=
(\partial/\partial z_\zs{m,j}) \tau_\zs{m,j}(z)
=\sqrt{h} \ov{e}_\zs{j}(\chi_\zs{\eta})\,,
$$
where 
\begin{equation}\label{sec:Lo.6}
\ov{e}_\zs{j}(f)=\int_{-1}^1\,e^2_\zs{j}(v)\,f(v)\,\d v\,.
\end{equation}
Now to obtain a lower bound for
the Bayes risk 
$\cE_\zs{n}(\wh{S}_\zs{n})$
 we make use of Theorem~\ref{Th.Tr.1} which implies that
\begin{equation}\label{sec:Lo.7}
\cE_\zs{n}(\wh{S}_\zs{n})\,\ge
\sum_{m=1}^{M_\zs{n}}\sum_{j=1}^{N_\zs{n}}\,
\frac{h \ov{e}^2_\zs{j}(\chi_\zs{\eta})}{F_\zs{m,j}+B_\zs{m,j}\,  
+t^{-2}_\zs{m,j}}\,,
\end{equation}
where $F_\zs{m,j}=\sum_{i=1}^n\,D^2_\zs{m,j}(x_\zs{i})\,
\E\,g^{-2}(x_\zs{i},S_\zs{\vartheta,n})$ and
$$
B_\zs{m,j}=
\frac{1}{2}\sum_{i=1}^n\,
\E\,
\frac{\wt{\L}^2_\zs{m,j}(x_\zs{i},S_\zs{\vartheta,n})}{g^4(x_\zs{i},S_\zs{\vartheta,n})}\,
$$
with
$\wt{\L}_\zs{m,j}(x,S)=\L_\zs{x,S} \left(D_\zs{m,j}\right)$.
In the Appendix we show that
\begin{equation}\label{sec:Lo.8}
\lim_\zs{n\to\infty}\,
\sup_\zs{1\le m\le M_\zs{n}}\,\sup_\zs{1\le j\le N_\zs{n}}
\left|
\frac{1}{nh}\,
F_\zs{m,j}\,-\,\ov{e}_\zs{j}(\chi_\zs{\eta}^2)\,
g^{-2}_\zs{0}(\wt{x}_m)
\right|=0
\end{equation}
and
\begin{equation}\label{sec:Lo.9}
\lim_\zs{n\to\infty}
\frac{1}{nh}
\sup_\zs{1\le m\le M_\zs{n}}\,\sup_\zs{1\le j\le N_\zs{n}}\,B_\zs{m,j}\,=\,0\,.
\end{equation}
This means that, for any $\nu>0$
and for sufficiently large $n$,
$$
\sup_\zs{1\le m\le M_\zs{n}}\,\sup_\zs{1\le j\le N_\zs{n}}
\frac{F_\zs{m,j}+B_\zs{m,j}+t^{-2}_\zs{m,j}}
{nh \ov{e}_\zs{j}(\chi_\zs{\eta}^2) g^{-2}_\zs{0}(\wt{x}_m)
+t^{-2}_\zs{m,j}} \le 1+\nu\,,
$$
where $\wt{x}_\zs{m}$ is defined in \eqref{Fa.4}. Therefore, if we  denote in \eqref{sec:Lo.7} 
$$
\kappa^2_\zs{m,j}=\,nh\,g^{-2}_\zs{0}(\wt{x}_m)\,t^2_\zs{m,j}
\quad\mbox{and}\quad
\psi_\zs{j}(\eta,y)
=\,
\frac{\ov{e}^2_\zs{j}(\chi_\zs{\eta})y}{\ov{e}_\zs{j}^2(\chi_\zs{\eta}^2)y+1}
$$
we obtain, for sufficiently large $n$,
$$
n^{\frac{2k}{2k+1}}\cE_\zs{n}(\wh{S}_\zs{n})\,\ge
\frac{n^{-\frac{1}{2k+1}}}{1+\nu}
\,
\sum_{m=1}^{M_\zs{n}}\,g^2_\zs{0}(\wt{x}_m)
\,
\sum_{j=1}^{N_\zs{n}}\,\psi_\zs{j}(\eta,\kappa^2_\zs{m,j})\,.
$$
In the Appendix we show that
\begin{equation}\label{sec:Lo.10}
\lim_\zs{\eta\to 0}\,\sup_\zs{N\ge 1}\sup_\zs{(y_\zs{1},\ldots,y_\zs{N})\in \bbr^{N}_\zs{+}}
\left|
\frac{\sum_{j=1}^{N}\,\psi_\zs{j}(\eta,y_\zs{j})}{\Psi_\zs{N}(y_\zs{1},\ldots,y_\zs{N})}
\,-\,1
\right|\,=\,0\,,
\end{equation}
where
$$
\Psi_\zs{N}(y_\zs{1},\ldots,y_\zs{N})=\sum^N_\zs{j=1}\,\frac{y_\zs{j}}{y_\zs{j}+1}\,.
$$
Therefore we can write that, for sufficiently large $n$,
\begin{equation}\label{sec:Lo.11}
n^{\frac{2k}{2k+1}}\cE_\zs{n}(\wh{S}_\zs{n})\,\ge
\frac{1-\nu}{1+\nu}\,n^{-\frac{1}{2k+1}}
\,
\sum_{m=1}^{M_\zs{n}}\,g^2_\zs{0}(\wt{x}_m)
\,\Psi_\zs{N_\zs{n}}(\kappa^2_\zs{m,1},\ldots,\kappa^2_\zs{m,N_\zs{n}})\,.
\end{equation}

Obviously, to obtain a "good" lower bound for the risk 
$\cE_\zs{n}(\wh{S}_\zs{n})$
 one needs to maximize the right-hand side of the inequality \eqref{sec:Lo.11}. Hence we
choose the coefficients $(\kappa^2_\zs{m,j})$ by maximization  the function
$\Psi_\zs{N}$, i.e.
$$
\max_\zs{y_\zs{1},\ldots,y_\zs{N}}\,\Psi_\zs{N}(y_\zs{1},\ldots,y_\zs{N})
\quad\mbox{subject to}\quad
\sum^N_\zs{j=1}y_\zs{j}j^{2k}\le R\,.
$$
The parameter $R>0$ will be chosen later to satisfy the condition $\A_\zs{3})$. 
By the Lagrange multipliers method it is easy to find that
the solution of this problem is given by
\begin{equation}\label{sec:Lo.12}
y^*_\zs{j}(R)=a^*(R)\,j^{-k}-1
\end{equation}
with
$$
a^*(R)=\frac{1}{\sum^N_\zs{i=1}\,i^{k}}
\left(R+\sum^N_\zs{i=1}i^{2k}\right)
\quad\mbox{and}\quad
1\le j\le N\,.
$$

To obtain a positive solution in \eqref{sec:Lo.12} we need to impose  the following condition
\begin{equation}\label{sec:Lo.13}
R> \,N^{k}\,\sum^N_\zs{i=1} i^{k}-\sum^N_\zs{i=1}i^{2k}\,.
\end{equation}
Moreover, from the condition $\A_\zs{3})$ we  obtain that
\begin{equation}\label{sec:Lo.14}
R\le \frac{2^{2k+1}(1-\ve)r\,
n\,h^{2k+1}_\zs{n}}{\pi^{2k}\wh{g}_\zs{0}}:=R^*_\zs{n}\,,
\end{equation}
where
$$
\wh{g}_\zs{0}=2h_\zs{n}\,\sum_{m=1}^{M_\zs{n}}\,g^2_\zs{0}(\wt{x}_m)\,.
$$
Note that by the condition $\H_\zs{4})$ the function 
$g_\zs{0}(\cdot)=g(\cdot,S_\zs{0})$ is continuous on 
the interval $[0,1]$, therefore 
\begin{equation}\label{sec:Lo.15}
\lim_\zs{n\to\infty}\wh{g}_\zs{0}=\int^1_\zs{0}g^2(x,S_\zs{0})\d x=\varsigma(S_\zs{0})
\quad\mbox{with}\quad
S_\zs{0}\equiv 0\,.
\end{equation}

Now we have to choose the sequence $(h_\zs{n})$. Note that if we put in \eqref{Fa.5}
\begin{equation}\label{sec:Lo.15-1}
t_\zs{m,j}=g_\zs{0}(\wt{x}_m)\sqrt{y^*_\zs{j}(R)}/\sqrt{nh_\zs{n}}
\quad\mbox{i.e.}\quad
\kappa�_\zs{m,j}\,=\,y^*_\zs{j}(R)\,,
\end{equation}
we can rewrite the inequality \eqref{sec:Lo.11} as
\begin{equation}\label{sec:Lo.16}
n^{\frac{2k}{2k+1}}\cE_\zs{n}(\wh{S}_\zs{n})\,\ge
\frac{(1-\nu)}{(1+\nu)}\,\frac{\wh{g}_\zs{0}\,\Psi^*_\zs{N_\zs{n}}(R)}{2h_\zs{n}}\, n^{-\frac{1}{2k+1}}\,,
\end{equation}
where
$$
\Psi^*_\zs{N}(R)=
N-\frac{\left(\sum^{N}_\zs{j=1}j^{k}\right)^2}{R+\sum^{N}_\zs{j=1}j^{2k}}\,.
$$
It is clear that
$$
k^2/(k+1)^2\le \liminf_\zs{N\to\infty}\inf_\zs{R>0}\Psi^*_\zs{N}(R)/N
\le \limsup_\zs{N\to\infty}\sup_\zs{R>0}\Psi^*_\zs{N}(R)/N\le 1\,.
$$
Therefore to obtain a positive finite asymptotic lower bound in \eqref{sec:Lo.16}
we have to take the parameter $h_\zs{n}$ as
\begin{equation}\label{sec:Lo.17}
h_\zs{n}=h_\zs{*}n^{-1/(2k+1)}N_\zs{n}
\end{equation}
with some positive  coefficient $h_\zs{*}$. Moreover,  the conditions 
\eqref{sec:Lo.13}-\eqref{sec:Lo.14} imply that,
for sufficiently large $n$,
$$
(1-\ve)r\,\frac{2^{2k+1}}{\pi^{2k}}
\,\frac{1}{\wh{g}_\zs{0}}\,h^{2k+1}_\zs{*}
>
\frac{1}{N^{k+1}_\zs{n}}\sum^{N_\zs{n}}_\zs{j=1}j^k-\frac{1}{N^{2k+1}_\zs{n}}
\sum^{N_\zs{n}}_\zs{j=1}j^{2k}\,.
$$
Moreover, taking into account that for sufficiently large $n$
$$
\wh{g}_\zs{0}\,
\frac{1}{N_\zs{n}}\sum^{N_\zs{n}}_\zs{j=1}
\left(
(j/N_\zs{n})^k
-
(j/N_\zs{n})^{2k}
\right)
<\,
\frac{(1+\ve)\,\varsigma(S_\zs{0})k}{(k+1)(2k+1)}\,,
$$
 we obtain the following condition on $h_\zs{*}$
\begin{equation}\label{sec:Lo.18}
h_\zs{*}\ge (\upsilon^*_\zs{\ve})^{1/(2k+1)}\,,
\end{equation}
where
$$
\upsilon^*_\zs{\ve}=\,\frac{(1+\ve)k}{c^*_\zs{\ve}(k+1)(2k+1)}
\quad\mbox{and}\quad
c^*_\zs{\ve}=\frac{2^{2k+1}(1-\ve)r}{\pi^{2k}\varsigma(S_\zs{0})}\,.
$$
To maximize the function $\Psi^*_\zs{N_\zs{n}}(R)$
on the right-hand side of the inequality \eqref{sec:Lo.16} we take  $R=R^*_\zs{n}$ defined in \eqref{sec:Lo.14}.
Therefore we obtain that
\begin{equation}\label{sec:Lo.19}
\liminf_\zs{n\to\infty}\,\inf_\zs{\wh{S}_\zs{n}}\,
n^{2k/(2k+1)}\cE_\zs{n}(\wh{S}_\zs{n})\,\ge
\,\varsigma(S_\zs{0})\,F(h_\zs{*})/2\,,
\end{equation}
where
$$
F(x)=\frac{1}{x}-\frac{2k+1}{(k+1)^2(c^*_\zs{\ve}(2k+1)x^{2k+2}+x)}\,.
$$
Furthermore, taking into account that
$$
F^{\prime}(x)=-
\frac{(c^*_\zs{\ve}(2k+1)(k+1)x^{2k+1}-k)^2}{(k+1)^2(c^*_\zs{\ve}(2k+1)x^{2k+2}+x)^2}\le 0
$$
we get
$$
\max_\zs{h_\zs{*}\ge (\upsilon^*_\zs{\ve})^{1/(2k+1)}}F(h_\zs{*})=
F((\upsilon^*_\zs{\ve})^{1/(2k+1)})=
\frac{(1+\ve^\prime)k}{k+1}
\,(\upsilon^*_\zs{\ve})^{-1/(2k+1)}\,,
$$
where
\begin{equation}\label{sec:Lo.19-1}
\ve^\prime=\frac{\ve}{2k+\ve k+1}\,.
\end{equation}

This means that to obtain in \eqref{sec:Lo.19} the maximal lower bound one has to take
in \eqref{sec:Lo.17}
\begin{equation}\label{sec:Lo.20}
h_\zs{*}=(\upsilon^*_\zs{\ve})^{1/(2k+1)}\,.
\end{equation}
It is important to note that if one defines the prior distribution $\mu_\zs{\vartheta}$
in the bayesian risk \eqref{sec:Lo.3} by formulas \eqref{Fa.5}, \eqref{sec:Lo.15-1}, \eqref{sec:Lo.17} and \eqref{sec:Lo.20}, then the bayesian risk would depend on a parameter
$0<\ve<1$, i.e.
$\cE_\zs{n}=\cE_\zs{\ve,n}$.

Therefore, the inequality \eqref{sec:Lo.19} implies that, for any $0<\ve<1$,
\begin{equation}\label{sec:Lo.21}
\liminf_\zs{n\to\infty}\,\inf_\zs{\wh{S}_\zs{n}}\,
n^{2k/(2k+1)}\cE_\zs{\ve,n}(\wh{S}_\zs{n})
\ge
\frac{(1+\ve^\prime)(1-\ve)^{1/(2k+1)}}{(1+\ve)^{1/(2k+1)}}\,\gamma_k(S_\zs{0})\,,
\end{equation}
where the function $\gamma_k(S_\zs{0})$ is defined in \eqref{sec:Ma.4} for $S_\zs{0}\equiv 0$.

Now to end the definition of the sequence of the random functions
$(S_\zs{\vartheta,n})$ defined by \eqref{Fa.4} and 
\eqref{Fa.5} one has to define the sequence 
$(N_\zs{n})$.
 Let us remember that we make use of the  sequence 
$(S_\zs{\vartheta,n})$ 
with the coefficients
$(t_\zs{m,j})$
constructed in \eqref{sec:Lo.15-1} for $R=R^*_\zs{n}$ given in 
\eqref{sec:Lo.14} and for the sequence $h_\zs{n}$ given by 
\eqref{sec:Lo.17} and \eqref{sec:Lo.20} for some fixed arbitrary $0<\ve<1$.

  We will  choose the sequence $(N_\zs{n})$ 
to satisfy the conditions $\A_\zs{1})$--$\A_\zs{4})$. One can take, for example, 
$N_\zs{n}=[\ln^{4} n]+1$. Then the condition $\A_\zs{1})$ is trivial. Moreover,
taking into account  that in
this case
$$
R^*_\zs{n}=\frac{2^{2k+1}(1-\ve)r}{\pi^{2k} \wh{g}_\zs{0}}\upsilon^*_\zs{\ve} N^{2k+1}_\zs{n}
=
\frac{\varsigma(S_\zs{0})}{\wh{g}_\zs{0}}
\frac{k}{(k+1)(2k+1)}\,N^{2k+1}_\zs{n}
$$
we find thanks to the convergence
\eqref{sec:Lo.15}
$$
\lim_\zs{n\to\infty}\,
\dfrac{R^*_\zs{n}+\sum^{N_\zs{n}}_\zs{j=1}j^{2k}}{N^{k}_\zs{n} \sum^{N_\zs{n}}_\zs{j=1}j^{k}}=1\,.
$$
Therefore, the solution \eqref{sec:Lo.12}, for sufficiently large $n$, 
satisfies the following inequality
$$
\max_\zs{1\le j\le N_\zs{n}}
y^*_\zs{j}(R^*_\zs{n})\,j^k
\le 2 N^k_\zs{n}\,.
$$
Now it is easy to see that  the condition 
 $\A_\zs{2})$ holds with $d_\zs{n}=\sqrt{N_\zs{n}}$ and 
the condition
$\A_\zs{4})$ holds for arbitrary 
$0<\epsilon_\zs{0}<1$.
As to the condition $\A_\zs{3})$, note that in view of the definition of $t_\zs{m,j}$ 
in \eqref{sec:Lo.15-1}
we get
\begin{align*}
\frac{1}{h_\zs{n}^{2k-1}}\,\sum^{M_\zs{n}}_\zs{m=1}\sum^{N_\zs{n}}_\zs{j=1}\,t^2_\zs{m,j}\,j^{2k}
&=\frac{1}{2 nh_\zs{n}^{2k+1}}\,\wh{g}_\zs{0}\,
\sum^{N_\zs{n}}_\zs{j=1}\,y^*_\zs{j}(R^*_\zs{n})\,j^{2k}\\
&
=\frac{R^*_\zs{n} \wh{g}_\zs{0}}{N^{2k+1}_\zs{n} 2\upsilon^*_\zs{\ve}}
=
(1-\epsilon)r
\left(\frac{2}{\pi}\right)^{2k}\,.
\end{align*}
Hence the condition $\A_\zs{3})$.

\medskip
\section{Estimation of non periodic function}\label{sec:Np}

Now we consider the estimation problem of the non periodic regression  function $S$ in the model
\eqref{I.1}. In this case we will estimate the function $S$ on any interior interval $[a,b]$
of $[0,1]$, i.e. for $0<a<b<1$.
 
It should be pointed out that at the boundary points $x=0$ and $x=1$, one must to make use of
kernel estimators (see Brua, 2007).

Let now $\chi$ be a infinitely differentiable $[0,1]\to\bbr_\zs{+}$ function such that 
$\chi(x)=1$ for $a\le x\le b$ and $\chi^{(k)}(0)=\chi^{(k)}(1)=0$ for all $k\ge 0$, for example,
$$
\chi(x)=\frac{1}{\eta}\,\int^{\infty}_\zs{-\infty}\,
V\left(\frac{x-z}{\eta}\right)\,
\Chi_\zs{[a',b']}(z)\,
\d z\,,
$$
where $V$ is some kernel function introduced in \eqref{Fa.1},
$$
a'=\frac{a}{2}\,,\quad b'=\frac{b}{2}+\frac{1}{2}
\quad\mbox{and}\quad
\eta=\frac{1}{4}\min(a\,,\,1-b)\,.
$$
Multiplying the equation \eqref{I.1} by the function $\chi(\cdot)$ and simulating the
i.i.d. $\cN(0,1)$ sequence $(\zeta_\zs{j})_\zs{1\le j\le n}$ one comes to the estimation
problem of the periodic regression function $\wt{S}(x)=S(x)\chi(x)$, i.e. 
$$
\wt{y}_\zs{j}=\wt{S}(x_\zs{j})+\wt{\sigma}_\zs{j}(S)\,\wt{\xi}_\zs{j}\,, 
$$
where 
$\wt{\sigma}_\zs{j}(S)=\sqrt{\sigma^2_\zs{j}(S)+\epsilon^2}$,
$$
\wt{\xi}_\zs{j}=\frac{\sigma_\zs{j}(S)}{\wt{\sigma}_\zs{j}(S)}\,\xi_\zs{j}+
\frac{\epsilon}{\wt{\sigma}_\zs{j}(S)}\,\zeta_\zs{j}\,.
$$
and $\epsilon>0$ is some sufficiently small parameter.

 It is easy to see that if the sequence $(\sigma_\zs{j}(S))_\zs{1\le j\le n}$ satisfies the conditions 
$\H_\zs{1})-\H_\zs{4})$, then the sequence $(\wt{\sigma}_\zs{j}(S))_\zs{1\le j\le n}$
satisfies these conditions as well with 
$$
\wt{\sigma}_\zs{j}(S)=\wt{g}(x_\zs{j},S)=\sqrt{g^2(x_\zs{j},S)\chi^2(x_\zs{j})+\epsilon^2}\,.
$$
\medskip
\section{Conclusions}\label{sec:Cn}
In conclusion, it should be noted that this paper completes the investigation of the
estimation problem of the nonparametric regression function for the heteroscedastic regression
model \eqref{I.1} in the case of quadratic risk. It is proved that the adaptive procedure
\eqref{sec:Ad.7} satisfies the non asymptotic oracle inequality and it is asymptotically
efficient for estimating a periodic regression function. Moreover, in Section~\ref{sec:Np}
we have explained how to apply the procedure to the case of non periodic function.
As far as we know, the procedure \eqref{sec:Ad.7} is unique for estimating the regression function at the model \eqref{I.1}.
 Let us remember once more the main steps of this investigation. The procedure \eqref{sec:Ad.7} combines the both principal aspects of nonparametric estimation: non asymptotic and asymptotic.
Non-asymptotic aspect is based on the selection model procedure with penalization (see e.g., Barron, Birg\'e and Massart, 1999, or  
Fourdrinier, Pergamenshchikov, 2007). Our selection model procedure differs from the commonly used one by a small coefficient in the penalty term going to zero that provides the
sharp non-asymptotic oracle inequality.
 Moreover, the commonly used selection model procedure is based on the least-squares estimators
whereas our procedure uses weighted least-squares estimators with the weights minimizing the
asymptotic quadratic risk that provides the asymptotic efficiency, as the final result. 
 From practical point of view, the procedure \eqref{sec:Ad.7} gives an acceptable accuracy
even for small samples as it is shown via simulations by  Galtchouk, Pergamenshchikov, 2008.

\medskip
\section{ Proofs}\label{sec:Pr}

\subsection{ Proof of Theorem~\ref{Th.U.1}}\label{subsec:Th.U.1}
To prove the theorem 
we will adapt to the  heteroscedastic case the corresponding proof
 from Nussbaum, 1985.

First, from \eqref{sec:Ad.2}  we obtain that, for any $p\in\cP_\zs{n}$,
\begin{equation}\label{sec:Up.1-1}
\E_\zs{S,p}\,\|\wt{S}-S\|_\zs{n}^2=
\sum_{j=1}^{n}\,(1\,-\,\wt{\lambda}_\zs{j})^2
\theta^2_\zs{j,n}
+
\frac{1}{n}
\sum_{j=1}^n\,\wt{\lambda}_\zs{j}^2\varsigma_\zs{j,n}(S)\,,
\end{equation}
where
 $$
\varsigma_\zs{j,n}(S)=\frac{1}{n}\,
\sum^n_\zs{l=1}\sigma^2_\zs{l}(S)\phi^2_\zs{j}(x_\zs{l})\,.
$$
Setting now
$\wt{\omega}=\omega(\wt{\alpha})$ with the function $\omega$ defined in \eqref{sec:Ad.4-1},
the index $\wt{\alpha}$ defined in \eqref{sec:Up.0-1},
$\wt{j}_0=[\wt{\omega} \ve_\zs{n}]$,
$\wt{j}_1=[\wt{\omega}/\ve_\zs{n}]$
and
$$
\varsigma_\zs{n}(S)=\frac{1}{n}\,
\sum^n_\zs{l=1}\sigma^2_\zs{l}(S)\,,
$$
we rewrite \eqref{sec:Up.1-1} as follows
\begin{align}\label{sec:Up.1-2}
\E_\zs{S,p}\,\|\wt{S}-S\|_\zs{n}^2&=\sum_{j=\wt{j}_0+1}^{\wt{j}_1-1}
(1\,-\,\wt{\lambda}_\zs{j})^2\theta^2_\zs{j,n}\\ \nonumber
&+\frac{\varsigma_\zs{n}(S)}{n}\,
\sum_{j=1}^n\,\wt{\lambda}_\zs{j}^2+\wt{\Delta}_\zs{1,n}+\wt{\Delta}_\zs{2,n}
\end{align}
with 
$$
\wt{\Delta}_\zs{1,n}\,=\,\sum_{j=\wt{j}_1}^{n}\,\theta^2_\zs{j,n}
\quad
\mbox{and}
\quad
\wt{\Delta}_\zs{2,n}\,=\frac{1}{n}\,\sum_{j=1}^n\,\wt{\lambda}^2_\zs{j}
\left(\varsigma_\zs{j,n}(S)-\varsigma_\zs{n}(S)\right)\,.
$$
Note that we have decomposed the first term on the right-hand of \eqref{sec:Up.1-1} into the sum
$$
\sum_{j=\wt{j}_0+1}^{\wt{j}_1-1}
(1\,-\,\wt{\lambda}_\zs{j})^2\theta^2_\zs{j,n}\,+\,\wt{\Delta}_\zs{1,n}\,.
$$
This decomposition allows us to show that $\wt{\Delta}_\zs{1,n}$ is negligible and further to
approximate the first term by a similar term in which the coefficients $\theta_\zs{j,n}$
will be replaced by the Fourier coefficients $\theta_\zs{j}$ of the function $S$.

Taking into account the definition of $\omega$ in \eqref{sec:Ad.4-1}
we can bound $\wt{\omega}$ as
$$
\wt{\omega} \ge (A_\zs{k})^{\frac{1}{2k+1}}\,
(n \ve_\zs{n})^{\frac{1}{2k+1}}\,.
$$
Therefore, by
Lemma~\ref{sec:Le.A.1} we obtain 
$$
\lim_{n\to\infty}\sup_{S\in W_r^k}\,n^{\frac{2k}{2k+1}}\,\wt{\Delta}_\zs{1,n}=0\,.
$$
Let us consider now the next term $\wt{\Delta}_\zs{2,n}$. We have
\begin{align*}
|\wt{\Delta}_\zs{2,n}|
&=\frac{1}{n^2}\,\left|\sum^n_{i=1}\,\sigma^2_\zs{i}(S)\,
\sum_{j=1}^n\,\wt{\lambda}^2_\zs{j}\,
(\phi^2_\zs{j}(x_\zs{i})-1) \right|\\
&\le \frac{\sigma_\zs{*}}{n}
\sup_\zs{0\le x\le 1}
\left|\sum_{j=1}^n\,\wt{\lambda}_\zs{j}^2\,
(\phi^2_\zs{j}(x)-1)\right|\,.
\end{align*}
Now by Lemma~\ref{sec:Le.A.2} and the
definition \eqref{sec:Ad.4}
we obtain directly the same property for $\wt{\Delta}_\zs{2,n}$, i.e.
$$
\lim_{n\to\infty}\sup_{S\in W_r^k}\,n^{\frac{2k}{2k+1}}\,|\wt{\Delta}_\zs{2,n}|=0\,.
$$
Setting 
$$
\wh{\gamma}_\zs{k,n}(S)=n^{\frac{2k}{2k+1}}
\sum_{j=\wt{j}_0}^{\wt{j}_1-1}(1-\wt{\lambda}_\zs{j})^2\theta^2_\zs{j}
+
\varsigma_\zs{n}(S) n^{-\frac{1}{2k+1}}
\sum_{j=1}^n\,\wt{\lambda}_\zs{j}^2
$$
and applying the well-known inequality
$$
(a+b)^2\le (1+\delta)a^2+(1+1/\delta)b^2
$$
to the first term on the right-hand side of the inequality \eqref{sec:Up.1-2} we obtain that, for any  $\delta>0$
and for any $p\in\cP_\zs{n}$,
\begin{align}\nonumber
\E_\zs{S,p}\,\|\wt{S}-S\|_\zs{n}^2 &\le(1+\delta)
\,\wh{\gamma}_\zs{k,n}(S)\,n^{-2k/(2k+1)}\\[2mm] \label{sec:Up.1-3}
&+\wt{\Delta}_\zs{1,n}+
\wt{\Delta}_\zs{2,n}+(1+1/\delta)\,
\wt{\Delta}_\zs{3,n}\,,
\end{align}
where 
$$
\wt{\Delta}_\zs{3,n}=\sum_{j=\wt{j}_0+1}^{\wt{j}_1-1}(\theta_\zs{j,n}\,-\,\theta_\zs{j})^2\,.
$$
Taking into account 
that $k\ge 1$ and that
$$
\wt{j}_\zs{1}\le 
\ov{\omega}\,\ve^{-1}_\zs{n}
+(A_\zs{k})^{\frac{1}{2k+1}}\,
n^{\frac{1}{2k+1}} 
(\ve_\zs{n})^{-(2k+2)/(2k+1)}\,,
$$
we can show through Lemma~\ref{sec:Le.A.3} 
that
$$
\lim_{n\to\infty}\sup_{S\in W_r^k}\,n^{\frac{2k}{2k+1}}\,\wt{\Delta}_\zs{3,n}=0\,.
$$
Therefore, the inequality  \eqref{sec:Up.1-3}
 yields
$$
\limsup_{n\to\infty} n^{\frac{2k}{2k+1}}
\sup_{S\in W^k_r}
\cR_\zs{n}(\wt{S},S)/\gamma_k(S)\le 
\limsup_{n\to\infty}
\sup_{S\in W^k_r}
\wh{\gamma}_\zs{k,n}(S)/\gamma_k(S)
$$
and to prove \eqref{sec:Up.1} it suffices to show that
\begin{equation}\label{sec:Up.3}
\limsup_{n\to\infty}
\sup_{S\in W^k_r}
\wh{\gamma}_\zs{k,n}(S)/\gamma_k(S)
\le 1\,.
\end{equation}
First, it should be noted 
that the definition \eqref{sec:Up.0}
and the
inequalities \eqref{sec:Co.6}-\eqref{sec:Co.6-1} imply directly
$$
\lim_\zs{n\to\infty}\sup_{S\in W^k_r}
\left|
\wt{t}_\zs{n}/\ov{r}(S)
-1
\right|=0\,.
$$
Moreover,
by the definition of $(\wt{\lambda}_\zs{j})_\zs{1\le j\le n}$ in 
\eqref{sec:Up.0-2},
for sufficiently large $n$, for which $\wt{t}_\zs{n}\ge \ov{r}(S)$
we find
\begin{align*}
\sup_\zs{j\ge 1}
\,n^{\frac{2k}{2k+1}}
\,
\frac{(1-\wt{\lambda}_\zs{j})^2}{(\pi j)^{2k}}
&= \pi^{-2k}
(A_\zs{k}\wt{t}_\zs{n})^{-2k/(2k+1)}
\le \pi^{-2k}
(A_\zs{k}\ov{r}(S))^{-2k/(2k+1)}\,.
\end{align*}
Therefore, by the definition of the coefficients $(a_\zs{j})_\zs{j \ge 1}$
in \eqref{sec:Co.4} 
$$
\limsup_{n\to\infty}
n^{\frac{2k}{2k+1}}
\sup_{S\in W^k_r}
\sup_{j\ge \wt{j}_0}\,\pi^{2k}
(A_\zs{k}\ov{r}(S))^{2k/(2k+1)}
(1-\wt{\lambda}_\zs{j})^2/a_\zs{j}
\le 1\,.
$$
Furthermore, in view of the definition \eqref{sec:Ad.4} we calculate directly
$$
\lim_{n\to\infty}
\sup_{S\in W^k_r}
\left|
\,
n^{-\frac{1}{2k+1}}\,\sum^n_\zs{j=1}\wt{\lambda}^2_\zs{j}
-
(A_\zs{k}\ov{r}(S))^{\frac{1}{2k+1}}
\int^1_\zs{0}(1-z^k)^2\d z
\right|=0\,.
$$
Now, the definition of $W^k_\zs{r}$ in \eqref{sec:Co.3}
and the condition \eqref{sec:Co.5} imply the inequality \eqref{sec:Up.3}. Hence Theorem~\ref{Th.U.1}.
\endproof

\medskip
\subsection{ Proof of Theorem~\ref{Th.M.3}}\label{subsec:Th.M.3}

In this section we prove Theorem~\ref{Th.M.3}.
 Lemma~\ref{sec:Le.U.1} implies that to prove the lower bounds \eqref{sec:Ma.6} and \eqref{sec:Ma.6-1}, it suffices
to show 
\begin{equation}\label{sec:Lo.1}
\liminf_{n\to\infty}\,\inf_{\wh{S}_\zs{n}}\,n^{\frac{2k}{2k+1}}\,
\cR_0(\wh{S}_\zs{n})\,\ge\,1\,,
\end{equation}
where
$$
\cR_0(\wh{S}_\zs{n})\,=\,
\sup_\zs{S\in W_r^k}\,\E_\zs{S,q}\,\|\wh{S}_\zs{n}-S\|^2/\gamma_k(S)\,.
$$

 For any estimator $\wh{S}_\zs{n}$, we denote by $\wh{S}^0_\zs{n}$ its projection onto $W_r^k$, i.e.\\
 $\wh{S}^0_\zs{n}=\hbox{\rm Pr}_\zs{W_r^k}(\wh{S}_\zs{n})$.
Since $W^k_\zs{r}$ is a convex set, we get 
$$
\|\wh{S}_\zs{n}-S\|^2\ge\|\wh{S}^0_\zs{n}-S\|^2\,.
$$
Now  we introduce the following set
\begin{equation}\label{Fa.11}
\Xi_\zs{n}=\{\max_\zs{1\le m\le M_\zs{n}}\,\max_\zs{1\le j\le N}\,\zeta^2_\zs{m,j}\le d_\zs{n}\}\,,
\end{equation}
where  $(\zeta_\zs{m,j})$ are i.i.d. $\cN(0,1)$ random variables 
from \eqref{Fa.5}
and the sequence $(d_\zs{n})_\zs{n\ge 1}$ is given in the condition $\A_\zs{2})$.
Therefore, we can write that
$$
\cR_0(\wh{S}_\zs{n})
\ge\int_{\{z:S_\zs{z,n}\in W^k_\zs{r}\}\cap\Xi_\zs{n}}
\,\frac{\E_\zs{S_\zs{z,n},q}\|\wh{S}^0_\zs{n}-S_\zs{z,n}\|^2}{\gamma_k(S_\zs{z,n})}\,\mu_\zs{\vartheta}(\d z)\,.
$$
Here the kernel function family $(S_\zs{z,n})$ is given in \eqref{Fa.4}
in which $N_\zs{n}=[\ln^4 n]+1$ and
 the parameter $h$ is defined in \eqref{sec:Lo.17} and \eqref{sec:Lo.20}; the measure
 $\mu_\zs{\vartheta}$ is defined in \eqref{sec:Lo.3}.
Moreover, note that on the set $\Xi$ the random function $S_\zs{\vartheta,n}$
is uniformly bounded, i.e.
\begin{equation}\label{Fa.18}
\|S_\zs{\vartheta,n}\|_\zs{\infty}=
\sup_\zs{0\le x\le 1}\,|S_\zs{\vartheta,n}(x)|\,
\le \sqrt{d_\zs{n}}\,t^*_\zs{n}\,,
\end{equation}
where the coefficient $t^*_\zs{n}$ is defined in \eqref{Fa.6}.
 
Thus, we estimate the risk $\cR_0(\wh{S}_\zs{n})$ from below as
$$
\cR_0(\wh{S}_\zs{n})
\ge
\frac{1}{\gamma^*_\zs{n}}\,
\int_{\{z:S_\zs{z,n}\in W^k_\zs{r}\}\cap\Xi_\zs{n}}\,
\E_\zs{S_\zs{z,n},q}\|\wh{S}^0_\zs{n}-S_\zs{z,n}\|^2\,\mu_\zs{\vartheta}(\d z)
$$
with 
\begin{equation}\label{sec:Lo.2}
\gamma^*_\zs{n}=\sup_\zs{\|S\|_\zs{\infty}\le \sqrt{d_\zs{n}}t^*_\zs{n}}\,\gamma_k(S)\,.
\end{equation}

By making use of the  Bayes risk \eqref{sec:Lo.3} with the prior distribution given by 
formulae 
\eqref{Fa.5}, \eqref{sec:Lo.15-1}, \eqref{sec:Lo.17}
and \eqref{sec:Lo.20} for any fixed parameter $0<\ve<1$
 we  rewrite  the lower bound for $\cR_0(\wh{S}_\zs{n})$ as 
\begin{equation}\label{sec:Lo.4}
\cR_0(\wh{S}_\zs{n})
\ge
\cE_\zs{\ve,n}(\wh{S}^0_\zs{n})/\gamma^*_\zs{n}-2\,
\Omega_\zs{n}/\gamma^*_\zs{n}
\end{equation}
with
$$
\Omega_\zs{n}=\E
(\Chi_\zs{\{S_\zs{\vartheta,n}\notin W^k_\zs{r}\}}\,+\,
\Chi_\zs{\Xi^c_\zs{n}})
(r+\|S_\zs{\vartheta,n}\|^2)\,.
$$

In Section~\ref{subsec:Br} we proved that the parameters in chosen prior distribution 
 satisfy the conditions
 $\A_\zs{1})$--$\A_\zs{4})$.
Therefore Propositions~\ref{sec:Pr.Fa.2}--\ref{sec:Pr.Fa.3} and the limit \eqref{Fa.12}
imply that, for any $\nu>0$,
$$
\lim_\zs{n\to\infty}\,n^{\nu}\,
\Omega_\zs{n}=0\,.
$$
Moreover, by the condition $\H_\zs{4})$ the sequence $\gamma^*_\zs{n}$ goes to $\gamma_k(S_\zs{0})$
as $n\to\infty$. Therefore, from this, \eqref{sec:Lo.21} and \eqref{sec:Lo.4} we get, for any $0<\ve<1$,
$$
\liminf_{n\to\infty}\,\inf_{\wh{S}_\zs{n}}\,n^{\frac{2k}{2k+1}}\,
\cR_0(\wh{S}_\zs{n})\,\ge\,\frac{(1+\ve^\prime)(1-\ve)^{\frac{1}{2k+1}}}{(1+\ve)^{\frac{1}{2k+1}}}\,.
$$
where $\ve^\prime$ is defined in \eqref{sec:Lo.19-1}.
Limiting here $\ve\to 0$ implies inequality \eqref{sec:Lo.1}. 
Hence Theorem~\ref{Th.M.3}.
\endproof


\renewcommand{\theequation}{A.\arabic{equation}}
\renewcommand{\thetheorem}{A.\arabic{theorem}}
\renewcommand{\thesubsection}{A.\arabic{subsection}}
\section{Appendix}\label{Se.A}
\setcounter{equation}{0}
\setcounter{theorem}{0}

\subsection{ Proof of Lemma~\ref{sec:Le.U.1}}\label{subsec:A.1}

First notice that, for any $S\in W_r^k$, one has
\begin{align*}
\|\wh{S}_\zs{n}-S\|_\zs{n}^2
=\|T_\zs{n}(\wh{S})-S\|^2+\Delta^*_\zs{1,n}+\Delta^*_\zs{2,n}\,,
\end{align*}
where
$$
\Delta^*_\zs{1,n}=2\sum_{j=1}^n\int_\zs{x_{j-1}}^{x_\zs{j}}(\wh{S}_\zs{n}(x_\zs{j})-S(x))(S(x)-S(x_\zs{j}))\d x
$$ 
and
$$
\Delta^*_\zs{2,n}=\sum_{j=1}^n\int_\zs{x_{j-1}}^{x_\zs{j}}\,(S(x)-S(x_\zs{j}))^2\d x\,.
$$

For any $0<\delta<1$, by  making use of the elementary inequality
 $$
2ab\le \delta a^2+\delta^{-1}b^2\,,
$$
 one gets
$$
\Delta^*_\zs{1,n}\,\le\, \delta\|T_\zs{n}(\wh{S})-S\|^2+\delta^{-1}\Delta^*_\zs{2,n}\,.
$$
Moreover,  for any $S\in W_r^k$ with $k\ge 1$,
by the Bounyakovskii-Cauchy-Schwarz inequality we obtain that
$$
\Delta^*_\zs{2,n}\le\frac{1}{n^2}\sum_{j=1}^n\int_\zs{x_{j-1}}^{x_\zs{j}}\,\dot{S}^2(t)\,\d t
=\frac{1}{n^2}\|\dot{S}\|^2\le\frac{r}{n^2}\,.
$$
Hence Lemma \ref{sec:Le.U.1}. 
\endproof

\subsection{ Proof of Theorem~\ref{Th.Tr.1}}\label{subsec:A.2}

For any $z=(z_\zs{1},\ldots,z_\zs{l})'\in\bbr^n$ we set 
$$
\varrho_\zs{i}(v,z)=\frac{1}{f(v,z)\Phi(z)}
\,\frac{\partial}{\partial z_\zs{i}}\left(f(v,z)\Phi(z)\right)\,.
$$
Note that due to the condition \eqref{sec:Co.6}, the density \eqref{subsec:Tr.2} is bounded, i.e.
$$
f(v,z)\le (2\pi g_\zs{*})^{-n/2}\,.
$$
So through \eqref{subsec:Tr.1-1} we obtain that
$$
\lim_\zs{|z_\zs{i}|\to\infty}
\tau(z)\,f(v,z)\varphi_\zs{i}(z_\zs{i})=0\,.
$$
Therefore, integrating by parts yields
\begin{align*}
\wt{\E}(\wh{\tau}_\zs{n}-\tau(\vartheta))\varrho_\zs{i}&=
\int_\zs{\bbr^{n+l}}
(\wh{\tau}_\zs{n}(v)-\tau(z))
\frac{\partial}{\partial z_\zs{i}}
\left(f(v,z)\Phi(z)\right)\d z\,\d v\\
&=
\int_\zs{\bbr^l}\left(\frac{\partial}{\partial\, z_\zs{i}}\tau(z)\right)
\Phi(z)\left(\int_\zs{\bbr^{n}} f(v,z)\d v \right)\d z=\ov{\tau}_\zs{i}\,.
\end{align*}
Now the Bounyakovskii-Cauchy-Schwarz inequality gives the following lower bound
$$
\wt{\E}(\wh{\tau}_\zs{n}-\tau(\vartheta))^2\ge
\ov{\tau}_\zs{i}^2/\wt{\E}\varrho_\zs{i}^2\,.
$$
To estimate the denominator in the last ratio, note that
$$
\varrho_\zs{i}(v,z)
=\wt{f}_\zs{i}(v,z)
+
\frac{\dot{\varphi}_\zs{i}(z_\zs{i})}{\varphi_\zs{i}(z_\zs{i})}
\quad\mbox{with}\quad
\wt{f}_\zs{i}(v,z)=
(\partial/\partial\,z_\zs{i})\ln f(v,z)\,.
$$
>From \eqref{subsec:Tr.1} it follows that
$$
\wt{f}_\zs{i}(v,z)=
\sum^{n}_\zs{j=1}(\xi^2_\zs{j}-1)\,
\frac{1}{2\sigma^2_\zs{j}(z)}
\frac{\partial}{\partial\,z_\zs{i}}\,\sigma^2_\zs{j}(z)
+
\sum^{n}_\zs{j=1}\,\xi_\zs{j}\,
\frac{S^{\prime}_\zs{i}(x_\zs{j})}{\sigma_\zs{j}(z)}\,.
$$
Moreover, the conditions $\H_\zs{2})$ and \eqref{subsec:Tr.2-1} imply
\begin{align*}
(\partial/\partial\,z_\zs{i})\,\sigma^2_\zs{j}(z)\,=\,
(\partial/\partial\,z_\zs{i})\,g^2(x_\zs{j},S_\zs{z})\,=\,
\wt{\L}_\zs{i}(x_\zs{j},
z)
\end{align*}
from which it follows
$$
\wt{\E}\,\left(
\wt{f}_\zs{i}(Y,\vartheta)
\right)^2\,
=F_\zs{i}+B_\zs{i}\,. 
$$
This implies inequality \eqref{subsec:Tr.3}. Hence Theorem~\ref{Th.Tr.1}.
\endproof

\medskip
\subsection{ Proof of Proposition~\ref{sec:Pr.Fa.1}}\label{subsec:A.3}

 First note that, for $0\le x\le 1$, we can represent the $l$th derivative as
\begin{equation}\label{Fa.8}
S_\zs{z,n}^{(l)}(x)=
\frac{1}{h^l}\sum_{m=1}^{M_\zs{n}}\,\sum^l_\zs{i=0}\,
\left(^{l}_{i}\right)\,\chi_\zs{\eta}^{(l-i)}(v_\zs{m}(x))
\,Q_\zs{i,m}(z,v_\zs{m}(x))\,,
\end{equation}
where 
$$
Q_\zs{i,m}(z,v)=
\sum_{j=1}^{N_\zs{n}}\,z_\zs{m,j}e^{(i)}_\zs{j}(v)\,.
$$
Therefore 
\begin{align*}
\|S^{(l)}_\zs{z,n}\|^2=
\frac{1}{h_\zs{n}^{2l-1}}\sum^{M_\zs{n}}_\zs{m=1}
\int^1_\zs{-1}\,
\left(
\sum^l_\zs{i=0}\,
\left(^{l}_{i}\right)\,\chi_\zs{\eta}^{(l-i)}(v)
\,Q_\zs{i,m}(z,v)
\right)^2\,\d v
\end{align*}
and by the Bounyakovskii-Cauchy-Schwarz inequality we obtain that
\begin{equation}\label{Fa.9}
\|S^{(l)}_\zs{z,n}\|^2\le
\frac{C^*(l,\eta)}{h_\zs{n}^{2l-1}}
\sum^l_\zs{i=0}\,
\ov{Q}_\zs{i}(z)
\end{equation}
with $C^*(l,\eta)=\max_\zs{-1\le v\le 1}\,\sum^l_\zs{i=0}\,
\left(\left(^{l}_{i}\right)\,\chi_\zs{\eta}^{(l-i)}(v)\right)^2$
and
$$
\ov{Q}_\zs{i}(z)=\sum^{M_\zs{n}}_\zs{m=1}
\int^1_\zs{-1}
\,Q^2_\zs{i,m}(z,v)
\,\d v
\,.
$$
Now we show that, for any $0\le i \le k-1$ and $\delta>0$,
\begin{equation}\label{Fa.10}
\lim_\zs{n\to\infty}
n^\nu\,
\P\left(
\ov{Q}_\zs{i}(\vartheta)>\delta h_\zs{n}^{2k-1}
\right)=0\,.
\end{equation}
To this end  note that
\begin{align*}
\int^1_\zs{-1}
\,Q^2_\zs{i,m}(\vartheta,v)
\,\d v&=\sum^{N_\zs{n}}_\zs{j=1}\vartheta^2_\zs{m,j}
\int^1_\zs{-1}(e^{(i)}_\zs{j}(v))^2\,\d v\\
&\le \left(\frac{\pi}{2}\right)^{2i}
\sum^{N_\zs{n}}_\zs{j=1}t^2_\zs{m,j}\,j^{2i}\zeta^2_\zs{m,j}
\,.
\end{align*}
Therefore, taking into account
the definition of the set $\Xi_\zs{n}$ in \eqref{Fa.11},
the functions $\ov{Q}_\zs{i}(\vartheta)$ with $0\le i\le k-1$ can be estimated on this set  as
$$
\ov{Q}_\zs{i}(\vartheta)\le\,
\left(\frac{\pi}{2}\right)^{2(k-1)}\,
d_\zs{n}
\sum^{M_\zs{n}}_\zs{m=1}
\sum^{N_\zs{n}}_\zs{j=1}t^2_\zs{m,j}\,
j^{2(k-1)}\,
$$
and by  \eqref{Fa.7} we get, for any $\delta>0$ and
 sufficiently large $n$,
$$
\P\left(
\ov{Q}_\zs{i}(\vartheta)
>\delta h_\zs{n}^{2k-1}\right)
\le
\,
\P\left(\Xi^c_\zs{n}\right)\,.
$$
Moreover, for sufficiently large $n$, 
$$
\P\left(\Xi^c_\zs{n}\right)\le \,M_\zs{n}\,N_\zs{n}\, e^{-d_\zs{n}/2}\,.
$$
Therefore,  the condition $\A_\zs{1})$  implies
\begin{equation}\label{Fa.12}
\limsup_\zs{n\to\infty}\,n^\nu\,
\P\left(\Xi^c_\zs{n}\right)=0\,,
\end{equation}
for any $\nu>0$. Hence Proposition~\ref{sec:Pr.Fa.1}. \endproof
\medskip
\subsection{ Proof of Proposition~\ref{sec:Pr.Fa.2}}\label{subsec:A.4}

First of all we prove that for $\epsilon$ from the condition  $\A_\zs{3})$
\begin{equation}\label{Fa.13}
\lim_\zs{n\to\infty}\,n^\nu\,
\P\left(\|S^{(k)}_\zs{\vartheta,n}\|>\sqrt{(1-\epsilon/4)r}\right)\,=0\,.
\end{equation}
Indeed, 
 putting in \eqref{Fa.8} $l=k$ we can represent the $k$th derivative of 
$S_\zs{z,n}$ as follows
\begin{equation}\label{Fa.14}
S_\zs{z,n}^{(k)}(x)=S'_\zs{z,n}(x)+
S''_\zs{z,n}(x)
\end{equation}
with
$$
S'_\zs{z,n}(x)=
\frac{1}{h^k}\sum_{m=1}^{M_\zs{n}}\,\sum^{k-1}_\zs{i=0}\,
\left(^{k}_{i}\right)\,\chi_\zs{\eta}^{(k-i)}(v_\zs{m}(x))
\,Q_\zs{i,m}(z,v_\zs{m}(x))
$$
and
$$
S''_\zs{z,n}(x)=
\frac{1}{h^k}\sum_{m=1}^{M_\zs{n}}\,\chi_\zs{\eta}(v_\zs{m}(x))
\,Q_\zs{k,m}(z,v_\zs{m}(x))\,.
$$
First, note that, we can estimate the norm of 
$S'_\zs{z,n}(x)$ by the same way as in the inequality \eqref{Fa.9}, i.e.
$$
\|S'_\zs{z,n}\|^2
\le
\frac{C^*(k,\eta)}{h_\zs{n}^{2k-1}}
\sum^{k-1}_\zs{i=0}
\ov{Q}_\zs{i}(z)\,.
$$
By making use of \eqref{Fa.10} we obtain 
that, for any $p>0$ and for any $\delta>0$,
\begin{equation}\label{Fa.15}
\lim_\zs{n\to\infty}\,n^\nu\,
\P\left(\|S'_\zs{z,n}\|>\delta\right)=0\,.
\end{equation}

Let us consider now the last term in \eqref{Fa.14}. Taking into account that
 $0\le \chi_\zs{\eta}(v)\le 1$ we get
\begin{align*}
\|S''_\zs{z,n}\|^2&=\frac{1}{h_\zs{n}^{2k-1}}
\sum^{M_\zs{n}}_\zs{m=1}
\int^1_\zs{-1}\,\chi_\zs{\eta}^2(v)Q^2_\zs{k,m}(z,v)\d v\\
&\le
\left(
\frac{\pi}{2}
\right)^{2k}
 \frac{1}{h_\zs{n}^{2k-1}}
\sum^{M_\zs{n}}_\zs{m=1}
\sum^{N_\zs{n}}_\zs{j=1}t^2_\zs{m,j}\,j^{2k}\,
\zeta^2_\zs{m,j}\,.
\end{align*}
Therefore from the condition $\A_\zs{3})$ we get for sufficiently large $n$
$$
\|S''_\zs{z,n}\|^2\le (1-\epsilon/2)r+\left(
\frac{\pi}{2}
\right)^{2k}
\sum^{M_\zs{n}}_\zs{m=1}\ov{\zeta}_\zs{m}
:=(1-\epsilon/2)r+\left(\frac{\pi}{2}\right)^{2k}\,Z_\zs{n}
$$
with
$$
\ov{\zeta}_\zs{m}=\frac{1}{h_\zs{n}^{2k-1}}
\sum^{N_\zs{n}}_\zs{j=1}t^2_\zs{m,j}\,
j^{2k}
\wt{\zeta}_\zs{m,j}
\quad\mbox{and}\quad
\wt{\zeta}_\zs{m,j}=\zeta^2_\zs{m,j}-1\,.
$$
We show that for any $\nu>0$ and for any $\delta>0$
\begin{equation}\label{Fa.16}
\lim_\zs{n\to\infty}\,n^\nu\,
\P\left(
\,|Z_\zs{n}|
>\delta\right)=0\,.
\end{equation}
Indeed, by the Chebychev inequality for any $\iota>0$
\begin{equation}\label{Fa.17}
\P\left(
\,|Z_\zs{n}\,|
>\delta\right)\le\,
\E\,( Z_\zs{n} )^{2\iota}/\delta^{2\iota}\,.
\end{equation}
Note now that
 according to  the Burkholder-Davis-Gundy inequality for any $\iota>1$ there
exists a constant $B^*(\iota)>0$ such that
$$
\E\left( Z_\zs{n}\right)^{2\iota}\le B^*(\iota)\,
\E\,
\left(
\sum^{M_\zs{n}}_\zs{m=1}\ov{\zeta}^2_\zs{m}
\right)^{\iota}\,.
$$
Moreover, by putting
$$
\wt{\zeta}_\zs{*}=\max_\zs{1\le m\le M_\zs{n}}\,
\max_\zs{1\le j\le N_\zs{n}}\,\wt{\zeta}^2_\zs{m,j}
$$
we can estimate the random variable $\ov{\zeta}_\zs{m}$ as 
$$
\ov{\zeta}^2_\zs{m}\,\le\,
\frac{N_\zs{n}}{h_\zs{n}^{4k-2}}
\sum^{N_\zs{n}}_\zs{j=1}t^4_\zs{m,j}\,j^{4k}\,
\wt{\zeta}_\zs{*}\,.
$$
Therefore, by the condition $\A_\zs{4})$, for sufficiently large $n$, 
\begin{align*}
\E\left( Z_\zs{n}\right)^{2\iota}&\le B^*(\iota)\,N^{\iota}_\zs{n}\,
h^{\iota\epsilon_\zs{0}}_\zs{n}\,\E\,\wt{\zeta}^{\iota}_\zs{*}\\
&\le 
B^*(\iota)\,
\E\,(\zeta^2-1)^{2\iota}\,
M_\zs{n}N^{\iota+1}_\zs{n}
h^{\iota\epsilon_\zs{0}}_\zs{n}\,,
\end{align*}
where $\zeta\sim\cN(0,1)$.
Now the  condition $\A_\zs{1})$ implies, for sufficiently large $n$,
$$
\E\left( Z_\zs{n}\right)^{2\iota}\le n^{-\delta_\zs{1}\,(\iota\epsilon_\zs{0}-2)}\,.
$$
Thus, choosing in  \eqref{Fa.17}
$$
\iota >\nu/(\epsilon_\zs{0}\delta_\zs{1})+2/\epsilon_\zs{0}
$$
 we obtain the limiting equality 
\eqref{Fa.16} which together with \eqref{Fa.14}-\eqref{Fa.15}
implies \eqref{Fa.13}. Now it is easy to deduce that 
Proposition~\ref{sec:Pr.Fa.1} yields Proposition~\ref{sec:Pr.Fa.2}.
\endproof

\medskip
\subsection{ Proof of Proposition~\ref{sec:Pr.Fa.3}}\label{subsec:A.5}

First of all, we recall that, due to the condition $\A_\zs{2})$,
$$
\lim_\zs{n\to\infty}\,\sum^{M_\zs{n}}_\zs{m=1}\sum^{N_\zs{n}}_\zs{j=1}t^2_\zs{m,j}
\le\lim_\zs{n\to\infty}\,\frac{d_\zs{n}}{h^{2k-1}_\zs{n}}\,
\sum^{M_\zs{n}}_\zs{m=1}\sum^{N_\zs{n}}_\zs{j=1}t^2_\zs{m,j}\,j^{2(k-1)}=0\,.
$$
Therefore, taking into account that
\begin{equation}\label{Fa.17-1}
\|S_\zs{\vartheta,n}\|^2\le 
h_\zs{n}\sum^{M_\zs{n}}_\zs{m=1}\sum^{N_\zs{n}}_\zs{j=1}t^2_\zs{m,j}\zeta^2_\zs{m,j}\,,
\end{equation}
we obtain, for sufficiently large $n$, 
$$
\E\,\|S_\zs{\vartheta,n}\|^2\,
\left(
\Chi_\zs{\{S_\zs{\vartheta,n}\notin W^{k}_\zs{r}\}}
+
\Chi_\zs{\Xi^c_\zs{n}}
\right)
\,
\le 
\max_\zs{m,j}\,
\E\,\zeta^2_\zs{m,j}
\left(
\Chi_\zs{\{S_\zs{\vartheta,n}\notin W^{k}_\zs{r}\}}
+
\Chi_\zs{\Xi^c_\zs{n}}
\right)
\,.
$$
Moreover, for any $1\le m\le M_\zs{n}$ and $1\le j\le N_\zs{n}$,
we estimate the last term as
\begin{align*}
\E\,\zeta^2_\zs{m,j}
\left(
\Chi_\zs{\{S_\zs{\vartheta,n}\notin W^{k}_\zs{r}\}}
+
\Chi_\zs{\Xi^c_\zs{n}}
\right)
&\le\,n\,\P(S_\zs{\vartheta,n}\notin W^{k}_\zs{r})\\[2mm]
&+
n\,\P(\Xi^c_\zs{n})
+
2\E\,\zeta^2\,\Chi_\zs{\{\zeta^2\ge n\}}\,,
\end{align*}
where $\zeta\sim \cN(0,1)$. By applying now Proposition~\ref{sec:Pr.Fa.2}
and the limit 
\eqref{Fa.12}
we come to Proposition~\ref{sec:Pr.Fa.3}.
\endproof

\medskip
\subsection{ Proof of Proposition~\ref{sec:Pr.Fa.4}}\label{subsec:A.6}

Taking into account the inequality \eqref{Fa.18}
and the condition $\H_\zs{1})$
we obtain
\begin{align*}
\E\,\left|g^{-2}(x,S_\zs{\vartheta,n})- g^{-2}_\zs{0}(x)\right|&\le 
\max_\zs{\|S\|_\zs{\infty}\le \sqrt{d_\zs{n}}\,t^*_\zs{n}}\,| g^{-2}(x,S)-g^{-2}_\zs{0}(x)|\\
&+(2/g_\zs{*})\,
\P\left(\Xi^c_\zs{n}\right)\,.
\end{align*}

Conditions $\A_\zs{2})$ and $\H_\zs{4})$ together with the limit relation 
\eqref{Fa.12} imply Proposition~\ref{sec:Pr.Fa.4}.
\endproof

\subsection{ Properties of the trigonometric basis}\label{subsec:A.7}

\begin{lemma}\label{sec:Le.A.1} For any function $S\in  W_r^k$,
\begin{equation}\label{A.0-3}
\sup_{n\ge 1}\sup_{1\le m\le n-1}\,m^{2k}\,
\left(\sum_{j=m+1}^{n}\,\theta_\zs{j,n}^2
\right)
\,\le\,\frac{4r}{\pi^{2(k-1)}}\,.
\end{equation}
\end{lemma}

\begin{lemma}\label{sec:Le.A.2}
For any $m\ge 0$,
\begin{equation}\label{A.0-2}
\sup_{N\ge 2}\quad\sup_{x\in [0,1]}N^{-m}\left|\sum_{l=2}^{N}\,l^m
\,
\left( \phi^2_\zs{l}(x)-1 \right)\,
\right|\le 2^m\,.
\end{equation}
\end{lemma}
Proofs of Lemma~\ref{sec:Le.A.1} and Lemma~\ref{sec:Le.A.2} are given in Galtchouk, Pergamenshchikov, 2007.

\begin{lemma}\label{sec:Le.A.3}
Let $\theta_\zs{j,n}$ and $\theta_\zs{j}$ be the Fourier coefficients defined in \eqref{sec:Ad.1}
and \eqref{sec:Co.3-1}, respectively.
 Then, for $1\le j\le n$
and $n\ge 2$,
\begin{equation}\label{A.0-1}
\sup_{S\in W^1_r}\,|\theta_\zs{j,n}-\theta_\zs{j}|\,\le\, 2\pi\,\sqrt{r}\,j/n\,.
\end{equation}
\end{lemma}
\noindent{\bf Proof.}
Indeed, we have
\begin{align*}
|\theta_\zs{j,n}-\theta_\zs{j}|&\,=\,
\left|\sum_{l=1}^n\int_{x_\zs{l-1}}^{x_\zs{l}}\,\left(S(x_l)\phi_\zs{j}(x_l)-S(x)\phi_\zs{j}(x)\right)\d x\right|\\
&\le\,n^{-1}\,\sum_{l=1}^n\int_{x_\zs{l-1}}^{x_\zs{l}}\,
\left(|\dot{S}(z)\phi_\zs{j}(z)|\,+ \,|S(z)\dot{\phi}_\zs{j}(z)|\right)\d z
\\
&=\,n^{-1}\,\int_{0}^{1}\,
\left(|\dot{S}(z)|\,|\phi_\zs{j}(z)|\,+ \,|S(z)|\,|\dot{\phi}_\zs{j}(z)|\right)\d z\,.
\end{align*}
By making use of the Bounyakovskii-Cauchy-Schwarz inequality we get
\begin{align*}
|\theta_\zs{j,n}-\theta_\zs{j}|&\le\,
n^{-1}\,\left(\|\dot{S}\|\,\|\phi\|\,+\,\|\dot{\phi}\|\,\|S\|\right)\\
&\le\,
n^{-1}\,\left(\|\dot{S}\|\,+\,\pi\,j\,\,\|S\|\right)\,.
\end{align*}
The definition of class $W^1_r$  implies (\ref{A.0-1}).
Hence Lemma~\ref{sec:Le.A.1}.
\endproof

\subsection{ Proofs of \eqref{sec:Lo.8} and \eqref{sec:Lo.9}}\label{subsec:A.8}

First of all, note that Proposition~\ref{sec:Pr.Fa.4},
the condition \eqref{sec:Co.6} and the condition $\H_\zs{4})$ imply that
\begin{equation}\label{A.3-1}
\lim_\zs{n\to\infty}
\max_\zs{1\le m\le M_\zs{n}}\sup_\zs{0\le x\le 1}
\Chi_\zs{\{|v_m(x)|\le 1\}}
\,\E\left|g^{-2}(x,S_\zs{\vartheta,n})\,-\,g^{-2}_\zs{0}(\wt{x}_m)
\right|\,=\,0\,.
\end{equation}
Let us show now that for any  continuously differentiable function
$f$ on $[-1,1]$
\begin{equation}\label{A.3-2}
\lim_\zs{n\to\infty}\,\sup_\zs{1\le m\le M_\zs{n}}\,
\left|
\frac{1}{nh}\,\sum^{n}_\zs{i=1}\,f(v_m(x_\zs{i}))\Chi_\zs{\{|v_m(x_\zs{i})|\le 1\}}
-
\int^{1}_\zs{-1}f(v)\d v
\right|=\,0\,.
\end{equation}
Indeed, setting 
$$
\Delta_\zs{n,m}=\frac{1}{nh}\,\sum^{n}_\zs{i=1}\,f(v_m(x_\zs{i}))\Chi_\zs{\{|v_m(x_\zs{i})|\le 1\}}
-
\int^{1}_\zs{-1}f(v)\d v\,,
$$
 we deduce
\begin{align*}
\left|
\Delta_\zs{n,m}
\right|&=
\left|
\frac{1}{nh}\,\sum^{i^*}_\zs{i=i_\zs{*}}\,f(v_m(x_\zs{i}))
-
\int^{1}_\zs{-1}f(v)\d v
\right|\\
&\le
\sum^{i^*}_\zs{i=i_\zs{*}}\,\int^{v_m(x_\zs{i})}_\zs{v_m(x_\zs{i-1})}|f(v_m(x_\zs{i}))-f(z)|\d z
+\max_\zs{|z|\le 1}|f(z)|(2-v^*+v_\zs{*})\,,
\end{align*}
where $i_\zs{*}=[n\wt{x}_\zs{m}-nh]+1$, $i_\zs{*}=[n\wt{x}_\zs{m}+nh]$, 
$$
v_\zs{*}=( [n\wt{x}_\zs{m}-nh]+1-n\wt{x}_\zs{m})/(nh)
\quad
\mbox{and}
\quad
v^*=([n\wt{x}_\zs{m}+nh]-n\wt{x}_\zs{m})/(nh)\,.
$$ 
Therefore, taking into account that the derivative of the function $f$ is bounded
on the interval $[-1,1]$ we obtain that
$$
\left|
\Delta_\zs{n,m}
\right|\le
\frac{3\max_\zs{|z|\le 1}|\dot{f}(z)|+
2\max_\zs{|z|\le 1}|f(z)|}{nh_\zs{n}}\,.
$$
Taking into account the conditions on the sequence $(h_\zs{n})_\zs{n\ge 1}$
given in $\A_\zs{1})$ we obtain limiting equality \eqref{A.3-2} which together
with \eqref{A.3-1} implies \eqref{sec:Lo.8}.

Now we study the behavior of $B_\zs{m,j}$. Due to the inequality
\eqref{sec:Co.8}
we estimate the Fr\'echet derivative as
\begin{equation}\label{A.3-3}
|\wt{\L}_\zs{m,j}(x,S_\zs{\vartheta,n})|\le\,C^*\,
\left(|S_\zs{\vartheta,n}(x)D_\zs{m,j}(x)|+|D_\zs{m,j}|_\zs{1}
+
\|S_\zs{\vartheta,n}\|\,\|D_\zs{m,j}\|
\right)\,.
\end{equation}
Consider now the fisrt term on the right-hand side of this inequality.
We have
\begin{align*}
\E\,(S_\zs{\vartheta,n}(x)D_\zs{m,j}(x))^2&=
\E\,\left(\sum^{N_\zs{n}}_\zs{l=1}\,\vartheta_\zs{m,l}\,e_\zs{l}(v_\zs{m}(x))\right)^2
e^2_\zs{j}(v_\zs{m}(x))I^4_\zs{\eta}(v_\zs{m}(x))\\
&\le 
\sum^{N_\zs{n}}_\zs{l=1}\,t^2_\zs{m,l}\,\Chi_\zs{\{|v_\zs{m}(x)|\le 1\}}
\le (t^*_\zs{n})^2\Chi_\zs{\{|v_\zs{m}(x)|\le 1\}}\,.
\end{align*}
We recall that the sequence $t^*_\zs{n}$ is defined in \eqref{Fa.6}.
Therefore, property \eqref{A.3-2} implies 
$$
\max_\zs{1\le m\le M_\zs{n}}\,\max_\zs{1\le j\le N_\zs{n}}\,
\frac{1}{nh}\,\sum^{n}_\zs{i=1}\,
\E(S_\zs{\vartheta,n}(x_\zs{i})D_\zs{m,j}(x_\zs{i}))^2
=\O((t^*_\zs{n})^2)\,.
$$
As to the second term on the right-hand side of \eqref{A.3-3}, we get
$$
|D_\zs{m,j}|_\zs{1}=\int^1_\zs{0}|e_\zs{j}(v_\zs{m}(x))\,
\chi_\zs{\eta}(v_\zs{m}(x))|\d x
=h\int^1_\zs{-1}| e_\zs{j}(v)\,
\chi_\zs{\eta}(v)|\d v\le 2h\,.
$$
Similarly, $\|D_\zs{m,j}\|^2\le h$ and,  by \eqref{Fa.17-1}  
$$
\E\|S_\zs{\vartheta,n}\|^2\le h\sum^{M_\zs{n}}_\zs{m=1}
\sum^{N_\zs{n}}_\zs{j=1}\,t^2_\zs{m,j}\le (t^*_\zs{n})^2\,.
$$
Therefore,
$$
\frac{1}{nh}\max_\zs{1\le m\le M_\zs{n}}\,\max_\zs{1\le j\le N_\zs{n}}
|B_\zs{m,j}|\,=\,\O((t^*_\zs{n})^2+h_\zs{n})
$$
and the condition $\A_\zs{1})$ implies \eqref{sec:Lo.9}. \endproof

\subsection{ Proof of \eqref{sec:Lo.10}}\label{Su.A.5}

Indeed, by the direct calculation it easy to see that, for any $N\ge 1$ and for any vector
$(y_\zs{1},\ldots,y_\zs{N})'\in\bbr^N_\zs{+}$,
\begin{align*}
\left|
\frac{\sum_{j=1}^{N}\,\psi_\zs{j}(\eta,y_\zs{j})}{\Psi_\zs{N}(y_\zs{1},\ldots,y_\zs{N})}
\,-\,1\right|
\le\,
\frac{
\max_\zs{j\ge 1}\left(
|\ov{e}^2_\zs{j}(\chi_\zs{\eta})-\ov{e}_\zs{j}(\chi_\zs{\eta}^2)|
+|\ov{e}^2_\zs{j}(\chi_\zs{\eta})-1|
\right)
}{\min_\zs{j\ge 1}\,\ov{e}_\zs{j}(\chi_\zs{\eta}^2)}\,,
\end{align*}
where the operator $\ov{e}_\zs{j}(f)$ is defined in
 \eqref{sec:Lo.6}.
Moreover, we remember that $\int^1_\zs{-1}e^2_\zs{j}(v)\d v=1$.
Therefore, taking into account the  property
\eqref{Fa.2} we obtain \eqref{sec:Lo.10}.
\endproof

\vspace{10mm}
{\bf Acknowledgements}

We are grateful for the comments and stimulating questions by two anonymous Referees
which have led to considerable improvements.


\begin{flushright}
\begin{tabular}{lcl}
   L.Galtchouk                       &$\quad$& S. Pergamenshchikov              \\
 Department of Mathematics           &$\quad$& Laboratoire de Math\'ematiques Raphael Salem,\\      
 Strasbourg University               &$\quad$& Avenue de l'Universit\'e, BP. 12,            \\
 7, rue Rene Descartes               &$\quad$&  Universit\'e de Rouen,                  \\
 67084, Strasbourg, France           &$\quad$&  F76801, Saint Etienne du Rouvray, Cedex France.\\
 e-mail: galtchou@math.u-strasbg.fr  &$\quad$& Serge.Pergamenchtchikov@univ-rouen.fr         \\
\end{tabular}
\end{flushright} 

\end{document}